\documentclass[10pt]{article}
\renewcommand{\baselinestretch}{1.1}

\usepackage[plainpages=false,pdfpagelabels,colorlinks=true,linkcolor=mypurple,urlcolor=mypurple,citecolor=mypurple]{hyperref}

\usepackage{tptstyle}

\usepackage[nameinlink,capitalize]{cleveref}

\usepackage{fullpage}
\usepackage[font={footnotesize,it}]{caption}

\title{Optimally taming biases in black-box models for efficient semiparametric estimation}
\input{authors.def}
\author{%
\authorOneFirst{} \authorOneLast{}$^{\authorOneAffil}$ ~~~~~~
\authorTwoFirst{} \authorTwoLast{}$^{\authorTwoAffil}$ ~~~~~~
\authorThreeFirst{} \authorThreeLast{}$^{\authorThreeAffil}$ ~~~~~~
\authorFourFirst{} \authorFourLast{}$^{\authorFourAffil}$
\\\\$^{1}$ \affilOne
\\$^{2}$ \affilTwo}
\date{June, 2026}

\begin{document}

\maketitle

\enlargethispage{3\baselineskip}
\begin{abstract}{
\renewcommand{\baselinestretch}{1.00}
Modern semiparametric estimation often relies on flexible ``black-box" machine learning methods to estimate nuisance functions, raising a fundamental question: how do nuisance estimation errors propagate into inference for low-dimensional target parameters? The dominant paradigm, exemplified by double/debiased machine learning (DML), yields error bounds in which nuisance estimation errors enter multiplicatively. While widely adopted, it remains unclear whether this multiplicative-rate dependence is optimal for black-box models. In this paper, we start by revisiting the partial linear model $Y = \mu_0(X) + T \cdot \beta_0 + \varepsilon$ under a structure-agnostic setting, where the nuisance function $\mu_0$ is estimated using a generic machine learning model, with approximation error $\delta^\appr_\mu$ and stochastic error  $\delta_\mu^\stoc$. We show that the standard DML rate is not optimal in the regime where the auxiliary function $\mathbb{E}[T|X=x]$ cannot be consistently estimated. We propose a new estimator for $\beta_0$ that achieves a sharper rate of $n^{-1/2}+ \delta^\appr_\mu + (\delta_\mu^\stoc)^2$ and establish a matching lower bound demonstrating its optimality. Our results reveal a new principle: the first-order stochastic error of nuisance estimation can be eliminated without imposing any additional assumptions. This also leads to a revised tuning strategy favoring under-smoothing, where $ \delta^\appr_\mu  \asymp (\delta_\mu^\stoc)^2$,  rather than the classical bias-variance trade-off with $ \delta^\appr_\mu  \asymp \delta_\mu^\stoc$. Under mild additional conditions, the estimator is asymptotically normal with minimal asymptotic variance. The proposed method extends to a broad class of semi-parametric linear functional estimation problems, including average treatment effect estimation. Our results imply that popular orthogonal score methods in semiparametric estimation with black-box nuisance learners can be substantially improved.

}
\end{abstract}
\nopagebreak
  \def\firstkw{T}%
  \def\kwd#1{\if T\firstkw\def\firstkw{F}\else, \fi#1\ignorespaces}%
%
%
%
\kwd{adversarial estimation}
\kwd{black-box model}
\kwd{causal inference}
\kwd{partial linear model}
\kwd{semiparametric linear functional}
\kwd{structure-agnostic optimality}
\unskip.%
\endgroup

\newpage
{ 
\maintextrefs
\def\nnhalf{n^{ -\frac{1}{2}}}
\def\nhalf{n^{\frac{1}{2}}}
\def\ninv{n^{-1}}

\section{Introduction}
Modern machine learning (ML) models, including the random forest \citep{breiman2001random}, deep learning models \citep{lecun2015deep}, and foundation models \citep{hollmann2025accurate}, excel at estimating complex functions in high-dimensional settings \citep{schmidt2020nonparametric, muller2022transformers, fan2024factor}. 
Their flexibility and adaptivity make them highly appealing in semi-parametric estimation, where a low-dimensional target parameter depends on potentially high-dimensional and complex nuisance functions. 

A canonical framework for integrating ML into semi-parametric inference in a \emph{structure-agnostic} manner is double/debiased machine learning (DML) \citep{chernozhukov2017double, chernozhukov2022debiased}. By utilizing orthogonal score functions, DML reduces sensitivity to nuisance estimation errors. It is methodologically scalable and theoretically grounded, ensuring that nuisance estimation errors enter the final error only via multiplicative remainder terms. Conversely, several structure-aware methods demonstrate that when the nuisance function belongs to specific classes such as H{\"o}lder \citep{robins2008higher} or sparse linear \citep{zhang2014confidence, bradic2022testability, wang2024debiasd} classes, specialized procedures can outperform the generic DML. This raises a fundamental question: 

\begin{quote}
{\it Do the generic DML already attain the best possible performance in the structure-agnostic setting, or can we do better without additional structures such as H{\"o}lder-smoothness or sparsity?}
\end{quote}

In this paper, we begin by studying this question in the partial linear model, a canonical semiparametric problem with rich structure-aware literature \citep{robinson1988root, donald1994series, carroll1997generalized, fan2005profile, robins2008higher, zhang2014confidence, van2014asymptotically, cai2017confidence, javanmard2017debiasing, bellec2022biasing}. We introduce a \emph{Structure-Agnostic Debiasing via Empirical balancing} (SADE) method that reduces the impact of black-box ML estimation errors on semiparametric estimation. In the regime where the auxiliary function $\pi_0(x) = \mathbb{E}[T|X=x]$ cannot be consistently estimated, we establish the \emph{structure-agnostic optimal} rate in estimating the linear coefficient. Our contributions are threefold. First, we propose a new estimator that achieves a sharper rate than the standard DML benchmark by reducing the impact of nuisance estimation errors on the target parameter. Second, we develop a minimax framework to characterize the best attainable target parameter estimation error under black-box procedures. Third, we extend the same SADE principle into a broader class of semi-parametric functional estimation problems. 

We instantiate these general results with neural networks (NNs), a leading example of black-box ML learners. Neural networks have been widely used for nonparametric function estimation, including multi-dimensional regression \citep{schmidt2020nonparametric, kohler2021rate, fan2024noise}, high-dimensional variable selection \citep{fan2024factor}, and causality pursuit \citep{gu2024causality}. Despite their growing use in semiparametric problems \citep{zhong2024neural, farrell2021deep}, their theoretical role in these settings remains poorly understood. Our analysis shows that the proposed bias-reduction method reduces the sensitivity of NN-based semiparametric estimation to
architectural choices such as width and depth, thereby enabling reliable performance in multi- and high-dimensional settings.

\subsection{Partial linear model under the structure-agnostic setting}
Let $\mathcal{D} = \mathcal{D}_1 \cup \mathcal{D}_2 = \{D_i\}_{i=1}^n \cup \{D_i\}_{i=n+1}^{2n}$ be $2n$  independent and identically distributed realizations of $D=(X,T,Y)$ generated from
\begin{equation}
\label{eq:intro:model}
\begin{alignedat}{2}
    Y &= \beta_0 \cdot T + \mu_0(X) + \varepsilon &&\qquad \text{with} \qquad \mathbb{E}[\varepsilon|X, T] = 0, \\
    T &= \pi_0(X) + u &&\qquad \text{with} \qquad \mathbb{E}[u|X] = 0.
\end{alignedat} 
\end{equation} 
Here $Y$ is the response, and $(X, T)$ are covariates. The second equation is only a decomposition and impose no additional model restrictions. Our goal is to estimate the regression coefficient $\beta_0 \in \mathbb{R}$. Since $\beta_0$ is the target parameter, $\mu_0(\cdot)$ and $\pi_0(\cdot)$ are {\it nuisance functions}.  Such a model is often used in the causal inference literature, where $T$ corresponds to the treatment variable.  We will call $\pi_0(X)$ an auxiliary function, as its information can be used to estimate $\beta_0$. We study this problem in a structure-agnostic setting: no structural assumptions, such as additivity or (sparse) linearity, are imposed on $\mu_0$ and $\pi_0$. Instead, we assume that $\mu_0$ and $\pi_0$ can be estimated by black-box ML classes $\mathcal{G}_{\mu}$ and $\mathcal{G}_{\pi}$, up to approximation error (i.e. misspecification bias) and stochastic error (i.e. standard error). The objective is to determine the best possible accuracy for estimating $\beta_0$ using these generic learners and their error budgets. 

Following the standard statistical learning notations \citep{bartlett2005local, vladimir2006local, massart2006risk}, we characterize the learning difficulty of a black-box model by two quantities. For a class $\mathcal{G}_{\mathsf{h}}$ used to estimate a target function $\mathsf{h}_0$, we define
$\delta^\appr_{\mathsf{h}} := \inf_{g\in \mathcal{G}_{\mathsf{h}}} \|g - \mathsf{h}_0\|_2$ 
as the approximation error, and $\delta^\stoc_{\mathsf{h}}$ as the stochastic error\footnote{The stochastic error depends on the sample size $n$; we omit this dependence for notational simplicity. We also use the sample size $n$ for defining the stochastic error rather than the total sample size $2n$, thereby leaving room for sample splitting. } which is formally defined as the critical radius of the local Rademacher complexity of the class $\mathcal{G}_{\mathsf{h}}$.  For example, for the linear function class with $p$ basis, $\delta_{\mathsf{h}}^\stoc \asymp \sqrt{p/n}$. Standard empirical process arguments imply that least squares over $\mathcal{G}_{\mathsf{h}}$ can estimate $\mathsf{h}_0$ at the error rate $O(\delta_\mathsf{h}^\appr + \delta_{\mathsf{h}}^\stoc)$. Thus, one typically has 
$$\|\hat{\pi} - \pi_0\|_2 \lesssim \delta^{\appr}_\pi + \delta^{\stoc}_\pi \qquad \text{for} \qquad \hat{\pi} = \argmin_{\pi\in \mathcal{G}_\pi} \ninv \sum_{i=1}^n [T_i - \pi(X_i)]^2.$$  

The central question is: given the black-box ML classes $\mathcal{G}_{\mu}$ and $\mathcal{G}_\pi$, together with their corresponding approximation error and stochastic error rates, what is the best possible accuracy for estimating $\beta_0$?

\subsection{The DML pipeline}
\label{sec:dml-plm}
Before addressing optimality, we first recall the standard structure-agnostic benchmark provided by DML. It estimates the nuisance functions using one split of the sample, and plug in the estimates into an orthogonal score on the other split. The DML pipeline first uses $\mathcal{D}_2$ to estimate $\pi_0$ and $\mu_0$ by least squares $\hat{\pi} \in \argmin_{\pi \in \mathcal{G}_\pi}  |\mathcal{D}_2|^{-1} \sum_{(X, T) \in \mathcal{D}_2} \left[\pi(X) - T\right]^2$ and estimate $(\beta_0, \mu_0)$ by joint least squares,
\begin{align*}
    \hat{\beta}, \hat{\mu} &\in \argmin_{\beta \in \mathbb{R}, \mu \in \mathcal{G}_\mu} |\mathcal{D}_2|^{-1} \sum_{(X, T, Y) \in \mathcal{D}_2} \left[Y - \beta \cdot T - \mu(X)\right]^2.
\end{align*} 
The target parameter $\theta_0=\beta_0$ is then estimated on $\mathcal{D}_1$ by the double robust  \citep{robins1994estimation} estimator   
\begin{align}\label{eq:dml-est}
    \hat{\theta}_{\mathtt{DML}} = \frac{\ninv \sum_{i=1}^n (Y_i - \hat{\mu}(X_i)) (T_i - \hat{\pi}(X_i))}{\ninv \sum_{i=1}^n T_i (T_i - \hat{\pi}(X_i))}. 
\end{align} 
The orthogonal score removes the leading first-order error in estimating the nuisance and yields
\begin{align}
    \label{eq:rate-dml}
    \begin{split}
        |\hat{\theta}_{\mathtt{DML}} - \theta_0| &\lesssim \nnhalf + \|\hat{\mu} - \mu_0\|_2 \cdot \|\hat{\pi} - \pi_0\|_2 \\
        &\lesssim \nnhalf + \left(\delta_\mu^\appr + \delta_\mu^\stoc\right) \cdot \left(\delta_\pi^\appr + \delta_\pi^\stoc\right),
    \end{split}
\end{align} 
where the second inequality substitutes the black-box estimation errors of $(\hat{\mu},\hat{\pi})$.  Thus, for fixed learners $(\mathcal{G}_\mu,\mathcal{G}_\pi)$, DML provides a structure-agnostic benchmark depending only on the approximation and stochastic error budgets of the two nuisance classes. When choosing candidate classes, \eqref{eq:rate-dml} suggests the familiar bias-variance trade-offs $\delta_\mu^\appr \asymp \delta_\mu^\stoc$ and $\delta_\pi^\appr \asymp \delta_\pi^\stoc$, matching the model selection for the best accuracy in estimating the nuisance functions themselves.

The estimator \eqref{eq:dml-est}, which uses  $\hat{\pi}$ to construct the debiasing direction, is one representative of a broader orthogonal score literature. Subsequent work has extended this paradigm in several directions, for example, by automating the debiasing step using fixed bases \citep{chernozhukov2022automatic} or generic ML models through adversarial estimation \citep{chernozhukov2026adversarial}. Although these approaches differ in how the orthogonal correction is implemented, they follow the same principle and yield the same multiplicative dependence on nuisance estimation errors in \eqref{eq:rate-dml}. 

The multiplicative rate is known to be sub-optimal in some structure-aware settings. For example, if $\mu_0$ and $\pi_0$ are sparse linear functions in $d$-dimensional space with sparsity level $s_\pi$ and $s_\mu$, DML yields the rate $n^{-1/2} + \sqrt{s_\pi \cdot s_\mu} \log(d)/n$ \citep{belloni2014inference, farrell2015robust}, whereas the structure-aware optimal rate is $n^{-1/2} + (s_\pi \land s_\mu) \log(d)/n$ \citep{bellec2022biasing}. Hence, DML can be strictly sub-optimal when the two nuisance difficulties are imbalanced. Such examples do not resolve the structure-agnostic question, since a faster rate exploits additional structure unavailable to a generic black-box procedure. They do suggest, however, that the generic multiplicative rate \eqref{eq:rate-dml} may not be the final answer in the structure-agnostic setting. This leads to our main question.
\begin{question}
\label{ques:main}
    Can we obtain a provably better error bound than \eqref{eq:rate-dml} in the structure-agnostic setting?
\end{question}

\subsection{Messages from prior structure-aware results}

To approach this question, we first look to structure-aware results for guidance. While they do not directly resolve the structure-agnostic problem, they reveal which components of the DML multiplicative rate are intrinsic and which may be reducible. Using the same notation $(\delta_{\mathsf{h}}^\appr,\delta_{\mathsf{h}}^\stoc)$, we revisit the classical linear and sparse linear results for the partial linear model. The key lesson is that approximation and stochastic error affect the target estimator asymmetrically: structure-aware methods improve over DML mainly by reducing how stochastic nuisance error propagates into the final estimation error \citep{gu2025open}. This is why we need to disentangle $\delta_{\mathsf{h}}^\appr$ and $\delta_{\mathsf{h}}^\stoc$ in our formalization, rather than working with the combined nuisance error $\|\hat{h} - h_0\|_2 \asymp \delta_{\mathsf{h}}^\appr + \delta_{\mathsf{h}}^\stoc$. 

We begin with the fixed linear class $\mathcal{G}^{p}_{\mathrm{lin}} = \{g(x) = \sum_{l=1}^p a_l \phi_l(x), (a_1,\ldots, a_p)\in \mathbb{R}^p\}$ with stochastic error $\sqrt{p/n}$. When $\mathcal{G}_\mu=\mathcal{G}_\pi=\mathcal{G}^{p}_{\mathrm{lin}}$, a natural estimator is the joint least-squares estimator $\hat{\theta}_{\mathtt{LS}} = \hat{\beta}$ obtained by regressing $Y$ on $T$ and the basis functions in $\mathcal{G}_\mu$. 
One can get the error rate
\begin{align}
\label{eq:error-linear}
    |\hat{\theta}_{\mathtt{LS}} - \theta_0| \lesssim \nnhalf + \delta^\appr_\mu \cdot \delta^\appr_\pi + \nnhalf \cdot \delta_\mu^\stoc \qquad \text{if} ~~ \delta_\mu^\stoc = o(1),
\end{align} 
where $\delta_\mu^\stoc = \sqrt{p/n}$; see, e.g., \cite{donald1994series} and Section 7.3.3 of \cite{fan2020statistical}. This rate eliminates cross terms involving approximation and stochastic errors, such as $\delta_\mu^\stoc \cdot \delta_\pi^\appr$. The reason is the projection orthogonality of the best $L_2$ approximation onto
$\mathcal{G}^{p}_{\mathrm{lin}}$: for any $L_2$ function $g_0$, and $\bar{g} = \argmin_{g\in \mathcal{G}^{p}_{\mathrm{lin}}} \|g-g_0\|_2$, we have
$\mathbb{E}\{g(X)[g_0-\bar g](X)\}=0$ for all $g\in\mathcal{G}^{p}_{\mathrm{lin}}$. In the H{\"o}lder smooth setting, this fixed-basis orthogonality is the source of the structure-aware improvement over generic DML.  

Now consider the sparse linear model $\mathcal{G}_{\mathrm{slin}}^{p, s} = \left\{g(x)=\sum_{l=1}^p a_l \phi_l(x): \|a\|_0 \le s\right\}$,  
where the dictionary $\phi(x) = (\phi_1(x), \ldots, \phi_p(x))^\top$ is well-conditioned in that $\mathbb{E}[\phi(X) \phi(X)^\top]$ has eigenvalues bounded awary from 0 and $\infty$. The stochastic error for this class is $\sqrt{s\log(p)/n}$. When $\mathcal{G}_\mu=\mathcal{G}_{\mathrm{slin}}^{p, s_\mu}$ and $\mathcal{G}_\pi = \mathcal{G}_{\mathrm{slin}}^{p, s_\pi}$ are built on the same dictionary, debiased-Lasso procedures \citep{zhang2014confidence, van2014asymptotically, javanmard14confidence} yield
\begin{align}
\label{eq:error-sparse-linear}
    |\hat{\theta}_{\mathtt{DL}} - \theta_0| \lesssim \nnhalf + \delta^\appr_\mu \cdot \delta^\appr_\pi + \delta^\appr_\mu \cdot (\delta^\stoc_\pi)^2 + \left[\delta^\stoc_\mu\right]^2 
\end{align} 
when $\delta_\pi^\stoc = o(1)$. In the complimentary setting where $\mu_0$ is dense but $\pi_0$ is approximately sparse,  \cite{bradic2022testability} obtains the rate $n^{-1/2} + \delta^\appr_\mu \cdot \delta^\appr_\pi + \left[\delta^\stoc_\pi + \delta^\appr_\pi\right]^2$ 
after accounting for misspecification. Combining the two estimators can recover the ``minimum of square error'', $(s_\mu\land s_\pi) \log(p)/n$, when the model is well-specified $\delta_\mu^\appr = \delta_\pi^\appr = 0$. 

These examples convey two messages. First, the approximation multiplicative $\delta^\appr_\mu \cdot \delta^\appr_\pi$ appears unavoidable; related lower bounds follow from \cite{balakrishnan2023fundamental, jin2024structure, jin2025hard, jin2025sharp} after translating their formalization in the present setting; see \cref{remark:bala}. Second, methods differ mainly in how they propagate the stochastic error to the final estimation error. This distinction is most transparent when the auxiliary function $\pi_0$ is misspecified by $\mathcal{G}_\pi$, i.e., $\delta_\pi^\appr = \Omega(1)$. 

As illustrated in \cref{fig:rate-compare}, with controlled approximation error $\delta_\mu^\appr$ and stochastic error $\delta_\mu^\stoc$ budgets, the deterioration in the stochastic error tracks the increasing adaptivity of the ML model classes: fixed linear sieves are least adaptive, sparse linear classes are more adaptive, and generic black-box classes are most adaptive. One might therefore expect a fully structure-agnostic method to pay an additional price beyond the worst structure-aware case, represented by the sparse linear models. The main result of this paper later shows that this expectation is false in the regime $\delta_\pi^\appr = \Omega(1)$. 

\begin{figure}
\begin{center}
\begin{tikzpicture}
    \draw[->,myorange] (-1, -0.5) -- (-1, -3.5);
    \node[myorange] at (-1, -1.5) {more};
    \node[myorange] at (-1, -2.) {adaptive};
    \node[myorange] at (-1, -2.5) {class};
    \node[anchor=west] at (0, 0) {\textbf{ML class $\mathcal{G}_\mu$}};
    \node[anchor=west] at (3, 0) {\textbf{Method}};
    \node[anchor=west] at (6, 0) {\textbf{Upper bound on $|\hat{\theta}-\theta_0|$}};
    \node[anchor=west] at (0, -1) {Linear};
    \node[anchor=west] at (0, -2) {Sparse linear};
    \node[anchor=west] at (0, -3) {Generic};
    \node[anchor=west] at (3, -1) {Joint LSE};
    \node[anchor=west] at (3, -2) {Debiased-Lasso};
    \node[anchor=west] at (3, -3) {DML};
    \node[anchor=west] at (6, -1) {$n^{-1/2} + \delta_\mu^\appr + \myred{n^{-1/2} \cdot \delta^\stoc_\mu}$};
    \node[anchor=west] at (6, -2) {$n^{-1/2} + \delta_\mu^\appr + \myred{[\delta^\stoc_\mu]^2}$};
    \node[anchor=west] at (6, -3) {$n^{-1/2} + \delta_\mu^\appr + \myred{\delta^\stoc_\mu}$};
    \draw[->,myred] (11.5, -0.5) -- (11.5, -3.5);
    \draw[myblue] (7.4, -0.5) rectangle (8.2, -3.5);
    \node[myblue] at (7.8, -3.8) {invariant};
    \node[myred] at (11.5, -1.5) {worse};
    \node[myred] at (11.5, -2) {dependency};
    \node[myred] at (11.5, -2.5) {on $\delta^\stoc_\mu$};
\end{tikzpicture}
\end{center}
\caption{An illustration of error rates for different methods when the auxiliary function $\pi_0$ is not estimable, or its estimator is not used.}
\label{fig:rate-compare}
\end{figure}

\subsection{New insights and main messages}

The comparison in \cref{fig:rate-compare} suggests a pessimistic possibility: as the nuisance class becomes more adaptive, the stochastic error may have a larger impact on the final estimator. In particular, one might expect a fully structure-agnostic procedure to inherit the linear dependence  $\delta_\mu^\stoc$ from the DML benchmark when $\pi_0$ cannot be consistently estimated. Our result shows that this pessimism is unnecessary. 

We construct a new estimator that does not explicitly estimate the auxiliary function $\pi_0$ and improves on the generic DML rate. In the regime where $\pi_0$ cannot be consistently estimated, the DML multiplicative term in \eqref{eq:rate-dml} is effectively of order 
$\delta_\mu^\appr+\delta_\mu^\stoc$. In contrast, our SADE estimator yields
\begin{align}
\label{eq:intro-sa-plm}
    |\hat{\theta} - \theta_0|
    \lesssim
    \nnhalf + \delta_\mu^\appr + (\delta_\mu^\stoc)^2 .
\end{align}
Thus, the answer to \cref{ques:main} is affirmative: even without imposing additional structure on the nuisance functions, one can improve over the DML multiplicative rate. In particular, the first-order stochastic error of the outcome learner can be eliminated.  

We also show that \eqref{eq:intro-sa-plm} is not further improvable under a structure-agnostic minimax framework that disentangles the effects of the approximation and stochastic errors. In the misspecified $\pi_0$ regime, we prove a matching lower bound, establishing the structure-agnostic optimality of \eqref{eq:intro-sa-plm}. Moreover, the same lower bound continues to hold even when $\mathcal{G}_\mu$ is restricted to a sparse linear class. Hence, when $\delta_\pi^\appr \asymp 1$, there is no gap between the generic structure-agnostic setting and the worst structure-aware sparse linear case.

The rate \eqref{eq:intro-sa-plm} also changes how the black-box learners should be tuned when $\pi_0$ cannot be consistently estimated. Consider a family of outcome classes $\{\mathcal{G}_{\mu,s}: s\in\mathcal{S}\}$ indexed by a complexity parameter $s$, where larger $s$ decreases approximation error but increases stochastic error. The DML rate \eqref{eq:rate-dml} suggests tuning as if the goal were to estimate $\mu_0$ itself, leading to the usual bias-variance trade-off $\delta_{\mu,s}^\appr \asymp \delta_{\mu,s}^\stoc$. In contrast, \eqref{eq:intro-sa-plm} shows that for optimal target estimation requires under-smoothing:
$$
\delta_{\mu,s}^\appr \asymp (\delta_{\mu,s}^\stoc)^2.
$$
We illustrate this principle for NN classes under which $\mu_0$ admits a hierarchical composition structure \citep{bauer2019deep, schmidt2020nonparametric}. For a fully connected ReLU network with width $N$ and depth $L$, $\delta_{\mu,s}^\appr \asymp (NL)^{-2\gamma^*}$ and $\delta_{\mu,s}^\stoc \asymp (NL)/ \sqrt{n}$, where $\gamma^*$ is the dimension-adjusted degree of smoothness  of the hardest component \citep{fan2024noise}.  The target-optimal choice is therefore $NL \asymp n^{1/(2+2 \gamma^*)}$ for the semiparametric estimation, rather than $NL \asymp n^{1/(2+4 \gamma^*)}$, which is optimal for estimating the nonparametric $\mu_0$. Consequently, $\hat\theta$ can be $\sqrt{n}$-consistent even when $\hat\mu$ converges at a nonparametric rate.  Indeed, if $n^{1/(4\gamma^*)} \ll NL \ll n^{1/4}$, which requires $\gamma^* \geq 1$, then $|\hat{\theta} - \theta_0| = O(n^{-1/2})$. In this regime, asymptotic normality with the semi-parametric efficiency bound is also attainable.

Finally, this debiasing mechanism applies well beyond the partial linear framework. Under standard regularity conditions, it eliminates the first-order stochastic error of the nuisance function for a broader class of linear functionals. Our findings therefore move past a single optimality case to reveal a general principle for semi-parametric estimation paired with black-box machine learning.

\subsection{Organization} 

This paper is organized as follows. We present our methods and upper-bound result for the partial linear model in \cref{sec:plm}: we introduce our method in \cref{sec:method-plm}, provide the rationale why our method can eliminate the first-order stochastic error, followed by presenting the upper bound in estimating $\theta_0 = \beta_0$ in \cref{sec:plm-est-error}. We further show how asymptotic normality can be attained under additional but mild conditions in \cref{sec:plm-normality}, and use the example of a neural network class to illustrate optimal convergence rate and the additional technical conditions in \cref{sec:ex1}. The structure-agnostic minimax risk is then introduced, accompanied by a matching result in \cref{sec:sa-plm-mu}. In \cref{sec:general}, we generalize our method and estimation error result to other linear functional estimation tasks. In \cref{sec:further}, we offer a full but non-matching result for the partial linear model when $\pi_0$ can also be estimated by some machine learning model $\mathcal{G}_\pi$, which serves as a benchmark and starting point for future studies. All the proofs are collected in the supplemental materials. 

\noindent{\bf Notations.} In this paper, we use $\const$ to denote a universal constant that is related to the boundedness from below and above conditions. Define $[n]=\{1,\ldots, n\}$. We let $a\lor b = \max \{a, b\}$ and $a\land b = \min\{a, b\}$. We use $a \lesssim b$, $b \gtrsim a$, or $a = O(b)$ if there exists some constant $C$ depending on $\const$ such that $a \le Cb$ for any $n \ge C$. Denote $a \asymp b$ if $a\lesssim b$ and $a \gtrsim b$. For a function $f: \mathcal{Z} \to \mathbb{R}$, we denote $\|f\|_\infty = \sup_{z \in \mathcal{Z}} |f(z)|$.

\subsection{Related works}

A substantial literature has established rates sharper than the generic DML benchmark in semiparametric estimation. In smooth settings, higher-order influence functions have been used to obtain faster rates when the nuisance functions are H{\"o}lder smooth \citep{robins2008higher, robins2016technical, liu2017semiparametric}. In sparse settings, a parallel line of work extends the ``minimum of squared errors'' rate to other functionals, including average treatment effect estimation, when the outcome and propensity score satisfy sparse linear or generalized sparse linear structure \citep{athey2018approximate, tan2020model, ning2020robust, wang2024debiasd, celentano2023challenges}. These improvements, however, rely on additional structure -- most notably H{\"o}lder smoothness or sparse linearity. There have also been attempts to sharpen the DML benchmark by adding further calibration steps on top of black-box nuisance estimates \citep{van2014targeted, bonvini2024doubly, van2024doubly}, but the existing analyses rely on assumptions on the first-stage estimates $(\hat{\mu}-\mu_0,\hat{\pi}-\pi_0,\hat{\mu},\hat{\pi})$, and it is unclear whether such conditions hold for generic black-box learners. \cite{jin2025hard} obtains a sharper result for the partial linear model with black-box learners based on extra fundamental assumptions on treatment noise $u$, while it remains unclear how their approach extends to general linear functional estimation. 

Our proposed method is connected to the literature on weights for balancing the covariate \citep{hellerstein1999imposing, zubizarreta2015stable, imai2014covariate, zhao2019covariate, athey2018approximate, hirshberg2021augmented, fan2022optimal}. These methods view the propensity score as a balancing score and construct weights to reduce bias by balancing covariates, or more generally, a pre-determined function class. Our estimator extends this balancing principle to generic black-box learners. Rather than parameterizing the weights through a low-complexity model, we adopt fully flexible empirical weights and choose them through an adversarial criterion that balances all plausible directions of the first-stage error induced by the ML class. This particular surgery is designed to replace the cross multiplicative term $\delta_\pi^\appr \cdot \delta_\mu^\stoc$ by $[\delta_\mu^\stoc]^2$ because adopting a low-complexity model will inevitably bring about an approximation error term like $\delta_\pi^\appr$. In particular, unlike earlier approaches that balance coordinate \citep{imai2014covariate, zhao2019covariate, athey2018approximate} or a convex class \citep{hirshberg2021augmented}, which is equivalent to (implicitly) regressing the propensity score (or the Riesz representer in general) over a function class \citep{bruns2025augmented} via minimax duality, our method has no intent to consistently estimating the auxiliary function $\pi_0$. Instead, it is developed only for balancing the entire difference class associated with the black-box learner that estimates $\mu_0$.

\section{SADE for partial linear model}
\label{sec:plm}

We now introduce SADE for a partial linear model in the regime where the auxiliary function $\pi_0$ cannot be consistently estimated. The estimator and analysis
therefore use only the outcome class $\mathcal{G}_\mu$; no model for $\pi_0$ is used.  We begin with a regularity condition and then define the approximation and stochastic error budgets for $\mathcal{G}_\mu$.
\begin{condition}
\label{cond:reg}
    Recall the partial linear model in \eqref{eq:intro:model}, there exists some constant $\const\ge 1$ such that the following holds:
    \begin{itemize}
    \item[(1)] All the components in the model are bounded: $|\varepsilon| \lor |u| \le \const$, $|\beta_0| \le \const$, $\|\mu_0\|_\infty \lor \|\pi_0\|_\infty \le \const$.
    \item[(2)] The machine learning model we used is uniformly bounded by $\const$, i.e., $\|g\|_\infty \le \const$ for any $g\in \mathcal{G}_\mu$. 
    \item[(3)] The treatment variable $T$ can not be fully explained by $X$ in that $\mathbb{E}[u^2] \ge 1/\const$.
    \end{itemize}
\end{condition}

We define the stochastic error of $\mathcal{G}_\mu$ through localized population Rademacher complexity. 
\begin{definition}[Localized population Rademacher complexity]
For a given radius $\delta>0$, function class $\mathcal{F}$, and distribution $\nu$ on $Z$, define
\begin{align*}
    \mathsf{R}_{n,\nu}(\delta;\mathcal{F}) = \mathbb{E}_{Z,R}\left[\sup_{h\in \mathcal{F}, \|h\|_{L_2(\nu)} \le \delta} \left|\ninv \sum_{i=1}^n R_i h(Z_i)\right|\right],
\end{align*} where $Z=(Z_1,\ldots, Z_n)$ is an i.i.d. sample from distribution $\nu$, and independent of $R=(R_1,\ldots, R_n)$  which are i.i.d. Rademacher variables taking values in $\{-1,+1\}$ with equal probability.
\end{definition} 

Denote $\nu_x$ as the distribution of $X$, we define the approximation error and stochastic error for the machine learning model $\mathcal{G}_\mu$ as
\begin{align}
\begin{split}
    \delta_{\mu}^\appr &:= \inf_{g\in \mathcal{G}_{\mu}} \|g - \mu_0\|_{L_2(\nu_x)} \qquad \qquad \text{and}  \\
    \delta_{\mu}^\stoc &:= \inf\left\{\delta_s>\sqrt{\frac{\log(n)}{n}}: \sup_{\delta \ge \delta_s} \frac{\mathsf{R}_{n,\nu_x}(\delta; \partial \mathcal{G}_{\mu})}{\const \cdot \delta} \le \delta_s\right\},
\end{split} \label{eq:errors-g}
\end{align}
where $\partial \mathcal{G}_\mu = \{g - \tilde{g}: g, \tilde{g} \in \mathcal{G}_\mu\}$ is the difference of the machine learning model $\mathcal{G}_\mu$.
Here we assume w.l.o.g. that $\delta_{\mu}^\stoc \ge \sqrt{\log(n)/n}$ for ease of presentation. One can replace the $\sqrt{\log(n)/n}$ in the inequality by $0$ at the cost of substituting $\delta_{\mu}^\stoc$ with $\delta_{\mu}^\stoc + \sqrt{\log(n)/n}$ in the remaining places where $\delta_{\mu}^\stoc$ appear. 

\subsection{SADE estimator}
\label{sec:method-plm}

We construct the SADE estimator by sample splitting.  On $\mathcal{D}_2$, we obtain an initial outcome estimate by  
\begin{align}
\label{eq:est1-step1}
    (\hat{\beta}, \hat{g}) \in \argmin_{|\beta| \le \const, g\in \mathcal{G}_\mu} \ninv \sum_{i=n+1}^{2n} \left[Y_i - \beta \cdot T_i - g(X_i) \right]^2.
\end{align} 
Standard empirical process arguments yield $\|\hat{g} - \mu_0\|_2 \lesssim \delta_{\mu}^\appr + \delta_{\mu}^\stoc$ with high probability. 
On $\mathcal D_1$, we construct empirical debiasing weights $\hat a=(\hat a_1,\ldots,\hat a_n)$ through the minimax program
\begin{align}
\label{eq:est1-step2}
\begin{split}
    \hat{a} \in \argmin_{\|a\|_\infty \le M} &\frac{\lambda}{n} \sum_{i=1}^n \max\{\hat{v}(X_i, T_i), 0\} \cdot a_i^2 \\
    &\qquad \qquad + \sup_{\beta \in \mathbb{R}, f\in \partial \mathcal{G}_\mu} \left| \ninv\sum_{i=1}^n (\beta \cdot T_i + f(X_i)) \cdot a_i - \beta\right| - \ninv \sum_{i=1}^n f^2(X_i),
\end{split}
\end{align} 
Here $\lambda$ and $M$ are tuning parameters to be determined, $\partial \mathcal{G}_\mu$ is the difference of the ML model defined in \eqref{eq:est1-step2}, $\hat{v}(\cdot)$ is another estimate of some function using $\mathcal{D}_2$ that help us minimize the asymptotic variance in estimating $\theta_0$ under heteroskedastic noise; one can let $\hat{v} \equiv 1$ when $\varepsilon$ is independent of $(X, T)$, and would choose a consistent estimate of $v_0(X, T) = \mathbb{E}[Y|X,T]$ in general. Maximizing over $\beta \in \mathbb{R}$ is equivalent to set the constraint $n^{-1} \sum_{i=1}^n T_i a_i = 1$.

The final SADE estimator is
\begin{align}
\label{eq:est1-step3}
\hat{\theta} = \ninv \sum_{i=1}^n \left(Y_i - \hat{g}(X_i) \right) \cdot \hat{a}_i \qquad \text{with} \qquad \hat{a} = (\hat{a}_1,\ldots, \hat{a}_n) \text{ in } \eqref{eq:est1-step2}.
\end{align}

\subsection{The rationale and estimation error}
\label{sec:plm-est-error}

To illustrate the rationale behind our estimator $\hat{\theta}$, we focus on a simplified case where $\varepsilon$ is independent of $(X,T)$ and $\hat v(X_i,T_i)\equiv 1$. In this case, the minimax program in \eqref{eq:est1-step2} behaves like the following constrained formulation obtained through the Lagrange multiplier method:
\begin{align}
\begin{split}
    \hat{a} \in &\argmin_{a \in \mathcal{A}} \ninv \sum_{i=1}^n a_i^2 \qquad \text{where} \\
    &~~ \mathcal{A} = \Bigg\{ a=(a_1,\ldots, a_n) \in [-M, M]^n: \ninv \sum_{i=1}^n a_i T_i = 1, \\
    &\qquad \qquad \sup_{f \in \partial \mathcal{G}_\mu} \left|\ninv \sum_{i=1}^n f(X_i) \cdot a_i \right|
    - \ninv \sum_{i=1}^n f^2(X_i) \le \tilde{C} (\delta^\stoc_\mu)^2\Bigg\}
\end{split}
\label{eq:plm-constrained}
\end{align}
for a sufficiently large constant $\tilde{C}_1 =\poly(\const)$. As to be shown below, the constraints \eqref{eq:plm-constrained} can be viewed as a generic black-box method aiming to remove the first-order stochastic error by using empirical debiasing weights targeted at $\partial \mathcal G_\mu$. It is an analogue of the balancing step in debiased Lasso \citep{zhang2014confidence, van2014asymptotically}. 

We now sketch why  \eqref{eq:plm-constrained}, combined with  \eqref{eq:est1-step3}, yields the rate $n^{-1/2} + \delta_\mu^\appr + (\delta_\mu^\stoc)^2$. 
First, the feasible set $\mathcal A$ is nonempty. Let $\bar{D} = \ninv \sum_{i=1}^n (T_i - \pi_0(X_i)) \cdot T_i$, we claim that the oracle weights $\bar{a} \in \mathbb{R}^n$ with $\bar{a}_i = (T_i - \pi_0(X_i)) / \bar{D} = u_i / \bar{D}$ satisfies $\bar{a} \in \mathcal{A}$ for large $M$. Since $\mathbb E[uT]=\mathbb E[u^2]>0$, $\bar D$ is bounded away from zero with high probability. Hence $\|\bar a\|_\infty\lesssim 1$ for large $M$, and  $\ninv\sum_{i=1}^n \bar a_iT_i=1$ by construction. It remains to prove the last balancing constraint. Substituting $\bar a_i = u_i/\bar D$, this reduces to  
\begin{align*}
\sup_{g,\tilde g\in\mathcal G_\mu}
\left\{\left|\ninv\sum_{i=1}^n (g-\tilde g)(X_i)\cdot u_i\right|-\ninv\sum_{i=1}^n (g-\tilde g)^2(X_i)\right\}\lesssim(\delta_\mu^\stoc)^2,
\end{align*} 
which follows from a multiplier empirical-process bound for the difference class $\partial\mathcal G_\mu$ given $u_i$ are i.i.d. zero-mean random variables conditioned on any fixed $X_1,\ldots, X_n$. 

Denote $\|f\|_n = [\ninv \sum_{i=1}^n f^2(X_i)]^{1/2}$. Since $\hat a = (\hat{a}(X_1, T_1),\ldots, \hat{a}(X_n, T_n))$ in \eqref{eq:plm-constrained} minimizes the empirical quadratic objective over $\mathcal A$ and $\bar a\in\mathcal A$, we also have with high probability $\ninv \sum_{i=1}^n \hat{a}_i^2 \le \ninv \sum_{i=1}^n \bar{a}_i^2 \lesssim 1$.
Now fix $\tilde g\in\mathcal G_\mu$ such that $\|\tilde g-\mu_0\|_2\le 2\delta_\mu^\appr$. Using the model \eqref{eq:intro:model} and the normalization $\ninv\sum_{i=1}^n \hat{a}_i\cdot T_i=1$ given $\hat{a} \in \mathcal{A}$, we decompose
\begin{align*}
    \hat{\theta} - \theta_0 &= \ninv \sum_{i=1}^n (\varepsilon_i + \beta_0 T_i + \mu_0(X_i) - \hat{g}(X_i)) \cdot \hat{a}_i - \beta_0 \\
    &= \underbrace{\ninv \sum_{i=1}^n \varepsilon_i \hat{a}_i}_{\mathsf{T}_1} + \underbrace{\ninv \sum_{i=1}^n \left(\mu_0 - \tilde{g}\right)(X_i) \cdot \hat{a}_i}_{\mathsf{T}_2} + \underbrace{\ninv \sum_{i=1}^n \left(\tilde{g}-\hat{g}\right)(X_i) \cdot \hat{a}_i}_{\mathsf{T}_3}. 
\end{align*} The first term is a sum of independent mean-zero variables conditioned on $\{(X_i,T_i)\}_{i=1}^n$, so $|\mathsf T_1|\lesssim n^{-1/2}$. The second term is controlled by Cauchy-Schwarz as $|\mathsf T_2|\le \sqrt{\ninv\sum_{i=1}^n \hat{a}_i^2} \cdot \|\mu_0-\tilde g\|_n\lesssim\delta_\mu^\appr$. 

The third term is where the balancing constraint is used. Since $\hat a\in\mathcal A$, the following instance-dependent error bounds holds
\begin{align}
    \forall g_1, g_2 \in \mathcal{G}_\mu \qquad \left|\ninv \sum_{i=1}^n (g_1 - g_2)(X_i) \cdot \hat{a}_i \right| \le \|g_1 - g_2\|_n^2 + \tilde{C} (\delta_\mu^\stoc)^2. 
\end{align} The above inequality is an instance-dependent error bound in that $g_1, g_2$ appears in both sides -- it makes sure all the potential first-order stochastic errors are eliminated in the empirical distribution $\mathbb{P}_n = \ninv \sum_{i=1}^n \delta_{(X_i, T_i)}$ by our weights $\hat{a}_1,\ldots, \hat{a}_n$. Applying this with $(g_1,g_2)=(\tilde g,\hat g)$ gives $|\mathsf{T}_3| \lesssim \|\hat{g} - \tilde{g}\|_n^2 + (\delta_\mu^\stoc)^2 \lesssim [\delta_\mu^\appr + \delta_\mu^\stoc]^2$. Since the error budgets are bounded, the squared approximation error $[\delta_\mu^\appr]^2$ is absorbed into $\delta_\mu^\appr$, thus the final rate can be concluded by combining the error bounds on $\mathsf{T}_1$--$\mathsf{T}_3$.

Formally, the following theorem provides an oracle-type inequality for the original estimator we constructed, under which $\hat{a}$ solves a full minimax objective \eqref{eq:est1-step2} with hyperparameter $(\lambda, M)$ instead of the constrained formulation in \eqref{eq:plm-constrained}. The optimal rate can be obtained by choosing hyperparameters $M \asymp 1$ and any $\lambda \lesssim 1/\sqrt{n}$.   We note here that the optimal estimation error rate does not require trading off $\lambda$, one can even choose $\lambda=0$. Tuning $\lambda$ is for minimal asymptotic variance and thus is required in \cref{thm:normal}.

\begin{theorem}
\label{thm:est1}
    Consider the partial linear model \eqref{eq:intro:model} satisfying \cref{cond:reg}. There exists a constant $\tilde{C} = \poly(\const)$ such that if $\|\hat{v}\|_\infty \le \const$ and we pick $M \ge 2\const^2$, then with probability at least $1-3e^{-t}-2\exp(-n(\delta_\mu^\stoc)^2)$, our estimator in \eqref{eq:est1-step3} satisfies
    \begin{align*} \tilde{C}^{-1} |\hat{\theta} - \theta_0| \le M \left(\sqrt{\frac{{t}}{{n}}} + \delta_\mu^\appr\right) + (\delta^\stoc_\mu)^2 + \lambda ,
    \end{align*}
  where $\delta_\mu^\appr$ and $\delta_\mu^\stoc$ are given by \eqref{eq:errors-g}.
\end{theorem}

\begin{remark}[The knowledge of the $\delta_\mu^\appr$ and $\delta_\mu^\stoc$] We remark here that our estimator $\hat{\theta}$ is also agnostic to the approximation error and the stochastic error budget $(\delta_\mu^\appr,\delta_\mu^\stoc)$ for the best performance for given $\mathcal{G}_\mu$. So does this when we want our estimate of $\hat{\theta}$ to satisfy the asymptotic normality as shown in \cref{thm:normal}. However, the knowledge of $\delta_\mu^\appr$ and $\delta_\mu^\stoc$ is required to get the best performance for a family of $\mathcal{G}_{\mu,s}$ with hyperparameter $s$ that trades off $\delta_\mu^\appr$ and $\delta_\mu^\stoc$.
\end{remark}

\subsection{Asymptotic normality}
\label{sec:plm-normality}

We now return to the general heteroskedastic case and show that the estimator can be made asymptotically normal with minimal asymptotic variance. For this purpose, we use another black-box model $\mathcal{G}_v$ to estimate the conditional variance function and use it to weight the minimax program in \eqref{eq:est1-step2}. Throughout this subsection, we will use the notations of function class $\mathcal{G}_{\mu,n}$ and $\mathcal{G}_{v,n}$ with dependency on $n$ because our result is a convergence-to-distribution result as $n\to \infty$. 

To be specific, we use the estimated $\hat{\beta}, \hat{g}$ and the dataset $\mathcal{D}_2$ to get an estimate of the conditional variance function $v_0(X,T) = \mathbb{E}[\varepsilon^2|X,T]$, 
\begin{align}
\label{eq:est1-v}
    \hat{v} = \argmin_{v \in \mathcal{G}_{v,n}} \ninv \sum_{i=n+1}^{2n} \left| \left(Y_i - \hat{\beta} T_i - \hat{g}(X_i) \right)^2 - v(X_i, T_i)\right|^2,
\end{align} and use the following estimate of the asymptotic variance when asymptotic normality is attained
\begin{align}
    \hat{\sigma} = \sqrt{\ninv \sum_{i=1}^n \hat{v}(X_i, T_i) \cdot \hat{a}_i^2}. 
\end{align}

We impose the following condition with respect to the estimation of the conditional variance function. This condition is relatively mild: the first part is a regularity condition. The second part requires a consistent estimation of the conditional variance function $v_0$, which is always possible for a fixed $v_0$ \citep{stone77consistent}. 

\begin{condition}\label{cond:v}
    Letting $\nu_{x,t}$ be the joint distribution of $(X,T)$, the following holds
    \begin{itemize}
        \item[(1)] $\|v_0\|_\infty \lor \|v_0^{-1}\|_\infty \le \const$.
        \item[(2)] $\delta_{v, n} = o(1)$ with   
        \begin{align}
            \delta_{v, n} := \inf_{g \in \mathcal{G}_{v,n}} \|g - v_0\|_{L_2(\nu_{x, t})} + \inf \left\{\delta_s > 0: \sup_{\delta \ge \delta_s} \frac{\mathsf{R}_{n,\nu_{x, t}}(\delta; \partial \mathcal{G}_{v, n})}{\const \cdot \delta} \le \delta_s\right\}.
        \end{align}
    \end{itemize}
\end{condition}

The asymptotic normality can be attained under a condition for $\mathcal{G}_{\mu, n}$ that is slightly stronger than $\delta_\mu^\appr + (\delta^{\stoc}_\mu)^2 = o(n^{-1/2})$ as presented in the following condition. \cref{cond:mu-normal} additionally requires that the approximation of $\mu_0$ in $\mathcal{G}_{\mu,n}$ be $o(1)$ in $L_\infty$ norm. 
\begin{condition} \label{cond:mu-normal}
    The stochastic error $\delta_{\mu, n}^\stoc$ for $\mathcal{G}_{\mu, n}$ is of order $o(n^{-1/4})$. Moreover, 
    \begin{align*}
        \inf_{g \in \mathcal{G}_{\mu,n}} \left[\frac{\|g - \mu_0\|_{L_2(\nu_x)}}{n^{-1/2}} \lor \|g - \mu_0\|_{L_\infty(\nu_x)}\right] = o(1). 
    \end{align*}
\end{condition}

\begin{theorem}
\label{thm:normal}
    Consider the partial linear model \eqref{eq:intro:model} satisfying \cref{cond:reg}. Assume further \cref{cond:v} and \cref{cond:mu-normal}  hold. If we pick $2\const^2 \le M = O(1)$, and $\lambda$ satisfying
    \begin{align*}
        (\delta^\stoc_{\mu, n})^2 \ll \lambda \ll n^{-1/2} 
    \end{align*}
    then our estimator \eqref{eq:est1-step3} with $\hat{v}(\cdot)$ in \eqref{eq:est1-v} satisfies
    \begin{align}
    \label{eq:anormality}
        \frac{\sqrt{n}(\hat{\theta} - \theta_0)}{\hat{\sigma}} \overset{d}{\to} \mathcal{N}(0, 1) \qquad \text{with} \qquad \hat{\sigma} = \sigma + \op(1).
    \end{align} where $\sigma = \left(\mathbb{E}[\varepsilon^2 u^2]\right)^{1/2} / \mathbb{E}[u^2]$ is the square root of the semi-parametric efficient asymptotic variance for the partial linear model. 
\end{theorem}

\subsection{Application to semiparametric estimation using neural networks}
\label{sec:ex1}

We illustrate Theorems \ref{thm:est1}--\ref{thm:normal} using NNs as black-box learners. We show that the SADE estimator attains the stated optimal rate and asymptotic normality when $\mu_0$ belongs to the hierarchical composition model (HCM), a function class to which NNs are adaptive \citep{schmidt2020nonparametric, kohler2021rate,fan2024factor}.

We use fully connected  ReLU NNs, where $\mathrm{ReLU}(\cdot) = \max\{0, \cdot\}$. For depth $L$ and width $N$, a \emph{deep ReLU network} has the form 
\begin{align}
\label{eq:nn-architecture}
    g(\cdot) = T_{L+1} \circ \overline{\mathrm{ReLU}}_L \circ T_L \circ \bar{\mathrm{ReLU}}_{L-1} \circ \cdots \circ T_2 \circ \overline{\mathrm{ReLU}}_1 \circ T_1(\cdot).
\end{align} 
Here $T_{l}(z) = W_l z + b_l: \mathbb{R}^{d_l} \to \mathbb{R}^{d_{l+1}}$ is a linear map, with weight matrix $W_l \in \mathbb{R}^{d_{l}\times d_{l-1}}$ and bias vector $b_{l} \in \mathbb{R}^{d_{l}}$. The layer dimensions satisfy $(d_0,d_1\ldots, d_L, d_{L+1}) = (d, N, \ldots, N, 1)$, and $\overline{\mathrm{ReLU}}_l$ applies the ReLU activation entrywise to vectors in $\mathbb{R}^{d_l}$. The equal-width assumption is used only for notational simplicity.

\begin{definition}[Deep ReLU network class]
    Define the family of deep ReLU networks taking $d$-dimensional vector as input with depth $L$, width $N$, and truncated by $B$ as 
        $$\mathcal{F}_{\mathtt{nn}}(d, L, N, B) = \{\tilde{g}(\cdot) = \mathrm{Tc}_B(g(\cdot)): g(\cdot) \text{ in } \eqref{eq:nn-architecture}\},$$ 
        where $\mathrm{Tc}_B: \mathbb{R} \to \mathbb{R}$ is the truncation operator defined as $\mathrm{Tc}_B(z) = {\min}\{|z|, B\} \cdot \mathrm{sign}(z)$.
\end{definition}

We next introduce the H{\"o}lder smooth function and HCM used to characterize the complexity of $\mu_0$.

\begin{definition}[$(\beta,C)$-smooth function]
    Let $\beta = r+s$, where $r\ge 0$ is an integer and $0<s\le 1$. A $d$-variate function $f$ is $(\beta, C)$-smooth if for every non-negative sequence $\alpha \in \mathbb{N}^d$ such that $\sum_{j=1}^d \alpha_j = r$, the partial derivative $\partial^{\alpha} f=(\partial f)/(\partial x_1^{\alpha_1}\cdots x_d^{\alpha_d})$ exists and satisfies $|\partial^{\alpha} f(x) - \partial^{\alpha} f(z)| \le C \|x-z\|_2^s$, for $C\ge 0$. Let $\mathcal{F}_{\mathtt{HS}}(d, \beta, C)$ denote this function class.
\end{definition}

\begin{definition}[Hierarchical composition model (HCM) $\mathcal{F}_{\mathtt{HCM}}(d, l, \mathcal{O}, C)$]
\label{hcm}
    We define function class of HCM $\mathcal{F}_{\mathtt{HCM}}(d, l, \mathcal{O}, C)$ \citep{kohler2021rate} with $l, d \in \mathbb{N}^+$, $C\in \mathbb{R}^+$, and $\mathcal{O}$, a subset of $[1,\infty) \times \mathbb{N}^+$, in a recursive way as follows. Let $\mathcal{F}_{\mathtt{HCM}}(d, 0,\mathcal{O}, C)=\{h(\cdot):h(x)=x_j, j\in [d]\}$, and for each $l\ge 1$, 
    \begin{align*}
    \mathcal{F}_{\mathtt{HCM}}(d, l,\mathcal{O}, C) = \big\{&h: \mathbb{R}^d \to \mathbb{R}: h(x) = g(f_1(x),...,f_t(x))\text{, where} \\
    &~~~~~ g\in \mathcal{F}_{\mathtt{HS}}(t, \beta, C) \text{ with } (\beta, t)\in \mathcal{O} \text{ and } f_i \in \mathcal{F}_{\mathtt{HCM}}(d, l-1,\mathcal{O}, C)\big\}.
\end{align*} 
\end{definition}

Thus, HCM consists of finite compositions of low-dimensional smooth functions. Its effective smoothness is $\gamma^\star=\min_{(\beta,t)\in\mathcal{O}}\beta/t,$
and the minimax optimal \(L_2\) estimation risk over this class is
$n^{-\gamma^\star/(2\gamma^\star+1)}$
\citep{kohler2021rate, fan2024factor}. For example, if
$f(x)=f_1(x_1)+f_2\{f_3(x_2,x_3),f_4(x_1,x_4,x_5)\},$
and all component functions have bounded second derivatives, then the hardest component is the three-variate \(f_4\), giving \(\gamma^\star=2/3\).

We can derive the result based on the following condition on the data-generating process.
\begin{condition} \label{cond:nn-plm} For the partial linear model in \eqref{eq:intro:model}, assume that for some constant $c_1 > 1$:
\begin{itemize}
\item[(1)] (Function complexity) $\mu_0 \in \mathcal{F}_{\mathtt{HCM}}(d, l, \mathcal{O}, C_h)$ with $\max_{(\beta, t) \in \mathcal{O}} (\beta\lor t) \lor l \lor C_h \lor d \le c_1$ and $\gamma_\mu^\star = \min_{(\beta, t)\in \mathcal{O}} (\beta/t)$. $v_0(x, t) = \mathrm{Var}[Y|X=x, T=t]$ is a $c_1$-Lipschitz function. 
\item[(2)] (Boundedness) All the components in the model are bounded: $\|X\|_\infty \le c_1$, $|\varepsilon| \lor |u| \le c_1$, $|\beta_0| \lor \|\mu_0\|_\infty \lor \|\pi_0\|_\infty \lor \|v_0\|_\infty \le c_1$. The conditional variance function $v_0$ is further bounded from below $\|v_0^{-1}\|_\infty \le c_1$.
\item[(3)] (Machine learning class and hyper-parameter $M$) We adopt NN classes for $\mathcal{G}_\mu=\mathcal{F}_{\mathtt{nn}}(d, L, N, B)$ and $\mathcal{G}_v=\mathcal{F}_{\mathtt{nn}}(d+1, L, N, B)$ with same hyper-parameters $L,N$ and $B=c_1$, and set $M=2c_1^2$. 
\end{itemize}
\end{condition}

\begin{proposition}
\label{prop:nn-plm}
Under \cref{cond:nn-plm}, the SADE estimator $\hat{\theta}$ in \eqref{eq:est1-step3} with NN hyperparameters $N, L$ satisfying $N \land L\ge \tilde{C}_1 \log(n)$ and $NL \asymp n^{\frac{1}{2(\gamma_\mu^\star+1)}} \log^{\frac{2\gamma^\star_\mu-1}{(\gamma^\star_\mu+1)}}(n)$ satisfies, with probability at least $1-n^{-10}$:
\begin{align}
\label{eq:rate-nn}
    |\hat{\theta} - \theta_0| \le \tilde{C}_1  \left[\lambda + \sqrt{\frac{\log(n)}{n}} + \left(\frac{\log^6(n)}{n}\right)^{\frac{\gamma_\mu^\star}{\gamma_\mu^\star+1}}\right]
\end{align} for some constant $\tilde{C}_1$ dependent on $c_1$. Moreover, as long as $\gamma_\mu^\star > 1$ and $(\log^6(n)/n)^{\frac{\gamma_\mu^\star}{\gamma_\mu^\star+1}} \ll \lambda \ll n^{-1/2}$, the asymptotic normality in \eqref{eq:anormality} holds.  
\end{proposition}

The required $\|\cdot\|_\infty$ approximation error rate condition follows from NN approximation results for HCM;  see Theorem 3.4 of \cite{fan2024noise}. Up to logarithmic factors, 
\begin{align*}
    \delta_{\mu}^\appr \asymp (NL)^{-2\gamma_\mu^\star} \qquad \text{and} \qquad \delta_\mu^\stoc \asymp \sqrt{\frac{(NL)^2}{n}}.
\end{align*} 
The oracle inequality \eqref{thm:est1} suggests the under-smoothing choice  $(NL)^2 \asymp n^{\frac{1}{(\gamma_\mu^\star+1)}}$, which makes  $\delta_\mu^\appr \asymp (\delta_\mu^\stoc)^2$. This differs from the choice of hyper-parameters $(NL)^2\asymp n^{\frac{1}{(2\gamma_\mu^\star+1)}}$, which balances $\delta_\mu^\appr \asymp \delta_\mu^\stoc$ for the optimal estimation of $\mu_0$. Consequently, the first step estimator is intentionally under-smoothed and attains only $\|\hat{g} - \mu_0\|_{L_2(\nu_x)} \asymp n^{-\frac{\gamma_\mu^\star}{2\gamma_\mu^\star+2}}$ which is slower than the optimal rate $ n^{-\frac{\gamma_\mu^\star}{2\gamma_\mu^\star+1}}$ for estimating $\mu_0$.

\begin{remark}[Minimax optimal rate for the HCM]
\label{remark:hcm}
When $\mu_0$ is $(\beta, C)$-smooth function and $\pi_0$ is not smooth, the minimax optimal rate in estimating $\theta_0=\beta_0$ is $n^{-1/2} + n^{-\frac{2\gamma_\mu}{2\gamma_\mu+1}}$ with $\gamma_\mu = \beta/d$ \citep{robins2008higher}. Here, the dependency on the dimension-adjusted-smoothness is sharper than the rate \eqref{eq:rate-nn}. It is then interesting to ask what the minimax optimal rate for estimating $\theta_0=\beta_0$ is when $\mu_0$ lies in HCM. Our lower bound in \cref{sec:sa-plm-mu} suggests that if one only use the facts $\delta_{\mu}^\appr \asymp (NL)^{-2\gamma_\mu^\star}$ and $\delta_\mu^\stoc \asymp \sqrt{\frac{(NL)^2}{n}}$, then \eqref{eq:rate-nn} is the best rate one can obtain. We conjecture that whether the minimax optimal rate in estimating $\mu_0$ can be improved to $n^{-1/2} + n^{-\frac{2\gamma_\mu^\star}{2\gamma_\mu^\star+1}}$ is related to whether one can find a series of fixed basis $\{\phi_l(x)\}_{l=1}^\infty$ such that
\begin{align*}
    \sup_{f_0 \in \mathcal{F}_{\mathtt{HCM}}(d, l, \mathcal{O}, C)} \inf_{\alpha = (\alpha_1,\ldots, \alpha_s)} \left\|f_0 - \sum_{j=1}^s \alpha_j \cdot \phi_j(x)\right\|_{L_2(\nu_x)} \lesssim s^{-\gamma_\mu^\star}.
\end{align*} 
While there is some work showing the wavelet basis cannot be adaptive to the HCM class \citep{schmidt2020nonparametric}, and some function classes cannot be approximated by any linear basis efficiently \citep{hayakawa2020minimax}, it is still open whether the HCM itself is a truly adaptive function set such that it cannot be learned by a series of linear basis efficiently. However, it is independent of the main question this paper studies; we leave it for future studies. 
\end{remark}


\begin{figure}
\begin{center}
\includegraphics[width=0.8\textwidth]{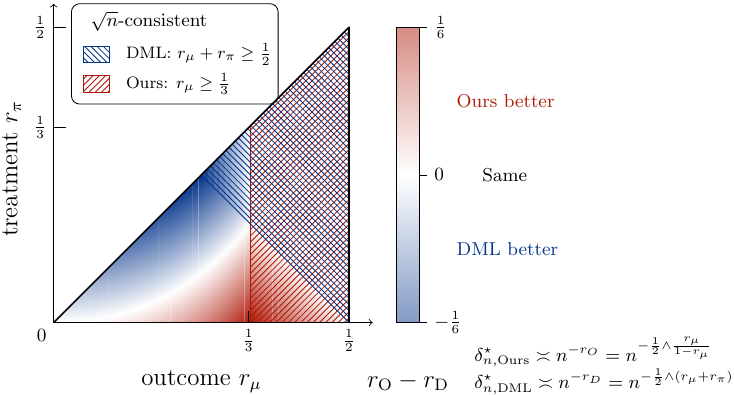}
\end{center}
\caption{Comparison of the best convergence rate $\delta^\star_n$ of $|\hat{\theta} - \theta_0|$ for our method (Ours) and double machine learning procedure (DML) when $\mu_0$ and $\pi_0$ both lie in HCM and the optimal rate in estimating $\mu_0$ and $\pi_0$ are $n^{-r_\mu}$ and $n^{-r_\pi}$ with $r_\mu, r_\pi \in [0,1/2]$. The color in $(r_\mu, r_\pi)$ represents the $r_O - r_D$, where $r_O$ (resp. $r_D$) is the $-\log_n(\text{best convergence rate }\delta_n^\star)$ for our methods (resp. DML) when the best rate in estimating $\mu_0$ and $\pi_0$ are $n^{-r_\mu}$ and $n^{-r_\pi}$ respectively. The red (resp. blue) slashed region represents the regime where our method (resp. DML) can attain $\sqrt{n}$ consistency. We can see that there is a considerable area (red) where our method outperforms DML, even though propensity score $\pi_0$ can also be learned at the rate of $n^{-r_\pi}$ with $r_\pi > 0$. }
\label{fig:rate}
\end{figure}

\begin{remark}[Comparison of our methods and DML when $\pi_0$ can also be estimated]
Although our method is structure-agnostically optimal only when $\pi_0$ cannot be consistently estimated, it offers a better estimation error rate than the DML in some regimes when $\pi_0$ can also be estimated well by some black-box learners. We illustrate this when both $\mu_0$ and $\pi_0$ are HCM and can be optimally estimated at the rate $n^{-r_\mu}$ and $n^{-r_\pi}$, respectively with $0\le r_\pi \le r_\mu \le 1/2$; see \cref{fig:rate}.
\end{remark}

\subsection{A numerical illustration}

We present a numerical example illustrating two implications of the theory in \cref{sec:ex1}: robustness to the difficulty of estimating the auxiliary function $\pi_0$, and the need for under-smoothing when estimating $\theta_0$. For $r>0$, define $f_r(v) = \{
\sum_{k=1}^{10} k^{-1/r} [\sin(k\pi \tanh(v)) - (-1)^k \cos(k\pi \tanh(v))]\} (\sum_{k=1}^{10}2k^{-2/r})^{-1/2}$. The normalization keeps the overall signal strength comparable across $r$. Larger $r$ yields slower coefficient decay and hence stronger high-frequency components, making $f_r$ harder to estimate.

We generate $X\sim \mathrm{Unif}([-1,1]^3)$, $(T,Y) \mid X$ from \eqref{eq:intro:model} with independent $\varepsilon, u \sim \mathrm{Unif}([-1, 1])$, $\theta_0 = \beta_0 \in[-1/2,1/2]$, $\mu_0(x)=f_1(w^\top x)$, $\pi_0(x)=f_r(w^\top x)$ with $w=\frac{1}{\sqrt{3}}(1,1,1)^\top$ and $r\in\{1,2,4,8\}$. Thus, $\mu_0$ is fixed, while $\pi_0$ becomes increasingly difficult to learn as $r$ increases. We set $n=|\mathcal{D}_1|=|\mathcal{D}_2|=1024$. The minimax objective \eqref{eq:est1-step2} can be solved by gradient descent ascent, which is widely used in the adversarial estimation literature \citep{goodfellow2020generative, arjovsky2017wasserstein, gu2024causality, chernozhukov2026adversarial}. Implementation details for our SADE method and DML are given in the supplement.

\begin{figure}
    \begin{center}
    \begin{tikzpicture}
        \node at (0, 0) {\includegraphics[height=4.1cm]{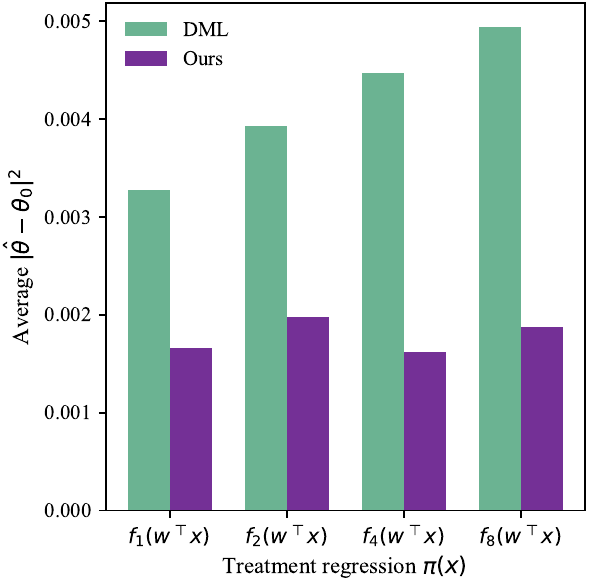}};
        \node at (6.9, 0) {\includegraphics[height=4.1cm]{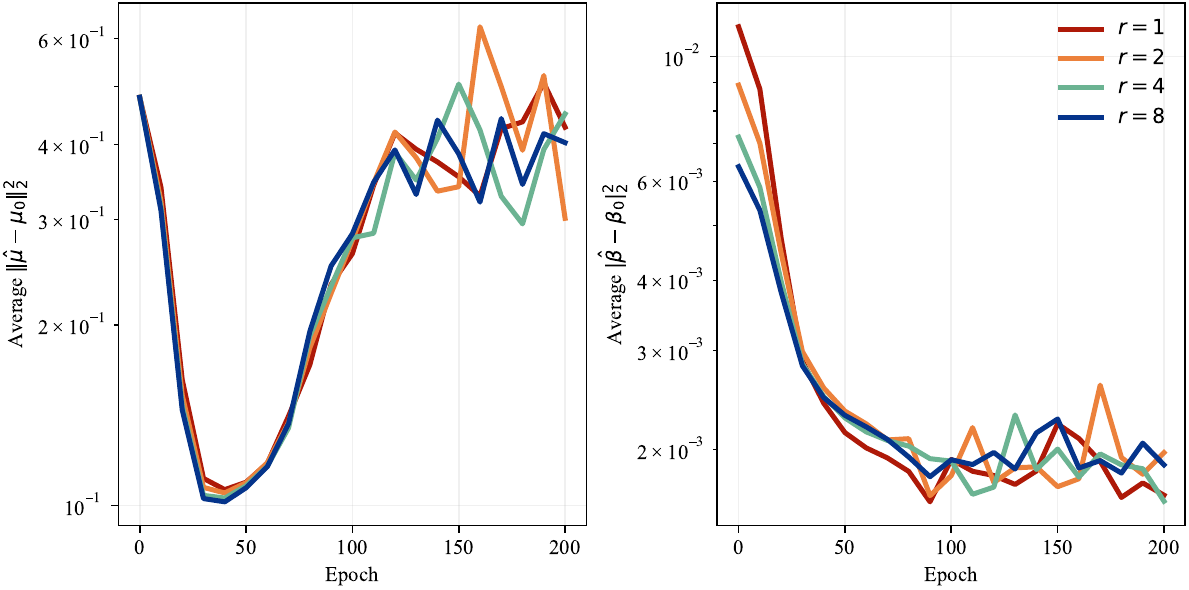}};
        \node at (0, -2.6) {(a)};
        \node at (6.9, -2.6) {(b)};
    \end{tikzpicture}
    \end{center}
    \caption{Numerical illustration of robustness to the auxiliary function difficulty and of under-smoothing. 
    Panel (a) depicts the average squared error $|\hat{\theta}-\theta_0|^2$ over 50 replications for SADE and DML as $r$ varies. Since $\mu_0$ is fixed and $\pi_0=f_r(w^\top x)$ becomes more complex as $r$ increases, the deterioration of DML reflects its sensitivity to the estimation error of $\pi_0$; whereas SADE remains stable across $r$. 
    Panel (b) reports the average nuisance estimation error $\|\hat{\mu}-\mu_0\|_2^2$ (left) and target estimation error $|\hat{\theta}-\theta_0|^2$ (right) over 50 replications as training epochs increase. The nuisance error is minimized at an earlier stopping time, whereas the target error is minimized after longer training, demonstrating that the tuning optimal for estimating $\theta_0$ under-smooths relative to the tuning optimal for estimating $\mu_0$ itself.}
    \label{fig:exp}
\end{figure}

\cref{fig:exp} (a) confirms the first phenomenon. Because only the complexity of $\pi_0$ changes, the increasing error of DML as $r$ grows reflects its multiplicative-remainder sensitivity to estimating $\pi_0$. In contrast, our SADE estimator does not require consistent estimation of $\pi_0$ in this regime; its empirical debiasing is targeted to the outcome learner, and the resulting estimation error remains nearly unchanged as $r$ varies.

\cref{fig:exp} (b) illustrates under-smoothing. The NN estimator of $\mu_0$ achieves its smallest $L_2$ error after roughly 40 epochs (left). The target estimation error $|\hat{\theta}-\theta_0|^2$ continues to decrease and is minimized after about 100 epochs (right). Hence, the best estimate of $\mu_0$ for optimal estimation of $\theta_0$ lies in the regime where the NN over-fits the data. This matches the theory: for estimating $\theta_0$, reducing approximation bias in $\hat{\mu}$ is more important than minimizing the full nuisance estimation error, because the first-order stochastic error is eliminated by debiasing.

\section{Structure-agnostic optimality}
\label{sec:sa-plm-mu}

In this section, we formalize structure-agnostic optimality and show that the rate in \cref{sec:plm} is sharp when the auxiliary function $\pi_0$ cannot be consistently estimated.  This is corresponds to restricting attention to procedures that do not exploit estimators of $\pi_0$. The quantity of interest is a minimax risk taken over all outcome classes that admit the same approximation and stochastic error budgets with respect to some black-box ML class.
For clarity, we state the lower bound under a simplified setting: $X\sim \nu_d :=\mathrm{Uniform}([0,1]^d)$, all nuisance functions and random variables are uniformly bounded by $3$, and $(\mathbb{E}[u^2])^{-1}\le 3$. These simplifications are not needed for the upper bound in \cref{sec:plm}; they are imposed only to simplify the matching lower bound. 

Let $\tfunc := \left\{f \in L_2(\nu_d): \|f\|_\infty \le 3\right\}$. For $\bar{\delta} = (\delta^\appr, \delta^\stoc) \in [0, 1]^2$, we define $\mathcal{H}_n(\bar{\delta})$ 
\begin{align}
\label{eq:h-n-delta}
    \mathcal{H}_n(\bar{\delta}) := \left\{\mathcal{F} \subseteq \tfunc: \exists \mathcal{G} \subseteq \tfunc, \sup_{f\in \mathcal{F}} \inf_{g\in \mathcal{G}} \|g - f\|_{L_2(\nu_d)} \le \delta^\appr,  \sup_{\delta\ge \delta^\stoc} \frac{\mathsf{R}_{n,\nu_d}(\delta; \partial \mathcal{G})}{3 \cdot \delta} \le \delta^\stoc\right\}.
\end{align} 
Thus, $\mathcal{H}_n(\bar{\delta})$ contains all target function classes that can be uniformly approximated to accuracy $\delta^\appr$ and learned with stochastic error $\delta^\stoc$ by some uniformly bounded black-box class using $n$ i.i.d. data.
For a given outcome function class $\mathcal{F}_{\mu}$, we define the associated family of distributions by
\begin{align}
\label{eq:dgp-plm-mu}
\begin{split}
    \mathcal{P}_\mu(\mathcal{F}_{\mu}) = \Big\{ &(X, T, Y) \sim \mathbb{P}_{\beta, v, \mu,\pi}~\text{where}~ \mu \in \mathcal{F}_{\mu}, \pi\in \tfunc, |\beta| \lor |v|\le 3, |T|\lor |Y| \le 3:\\
    &~~~~X\sim \nu_d, ~~ T = \pi(X) + u ~\text{with}~ \mathbb{E}[u|X]=0,\mathbb{E}[u^2|X] \equiv v \ge 1/3\\
    &~~~~Y = T \cdot \beta + \mu(X) + \varepsilon ~\text{with}~ \mathbb{E}[\varepsilon|X, T]=0\Big\}.
\end{split}
\end{align} 
The target parameter is $\theta(\mathbb{P}) = \beta$ for $\mathbb{P} = \mathbb{P}_{\beta, v, \mu, \pi} \in \mathcal{P}_\mu(\mathcal{F}_\mu)$. The restriction $\mathbb{E}[u^2| X]\equiv v$ is imposed only for the lower bound; it rules out additional learnable structure in the treatment variance. Throughout this section, $\pi_0$ is unrestricted over $\tfunc$, so it cannot be uniformly consistently estimated. Extensions to the setting where both $\mu_0$ and $\pi_0$ are learnable are given in \cref{sec:plm-full}.

For $\bar{\delta}_\mu=(\delta_\mu^\appr,\delta_\mu^\stoc)$, define the structure-agnostic minimax risk 
\begin{align}
\label{eq:sa-risk-mu}
    \mathfrak{m}(n, \bar{\delta}_\mu) := \sup_{\mathcal{F}_\mu \in \mathcal{H}_n(\bar{\delta}_\mu)} \underbrace{\inf_{\hat{\theta}} \sup_{\mathbb{P}\in \mathcal{P}_\mu(\mathcal{F}_\mu)} \mathbb{E}_{\mathbb{P}^{2n}}\left[\left|\hat{\theta} - \theta(\mathbb{P})\right|\right]}_{\text{standard minimax risk over given class $\mathcal{F}_\mu$}},
\end{align} 
where the randomness in $\mathbb{E}_{\mathbb{P}^{2n}}$ is $2n$ i.i.d. observations from $\mathbb{P}$, and $\hat{\theta}$ is the function of the $2n$ observations and potential black-box ML model $\mathcal{G}_\mu$ that is tailored to $\mathcal{F}_\mu$. Thus, $\mathfrak{m}(n,\bar{\delta}_\mu)$ is the worst standard minimax risk over all outcome classes sharing the same approximation and stochastic error budgets $\bar{\delta}_\mu$.

At a high level, the upper bound of the risk \eqref{eq:sa-risk-mu} can be proved by designing an algorithm having this error bound only using the conditions that $\mathcal{G}_\mu$ has stochastic error $\delta_\mu^\stoc$ and can approximate the target at the rate $\delta_\mu^\appr$, which is what we did in \cref{sec:plm}. On the other hand, the lower bound of the risk \eqref{eq:sa-risk-mu} can be established by finding a hard function class $\mathcal{F}_\mu$ within the approximation and stochastic errors budgets; the minimax optimal risk for the function class $\mathcal{F}_\mu$ can give a candidate lower bound for our structure-agnostic risk \eqref{eq:sa-risk-mu}. This idea will be made explicit as we present the matching result in the following theorem. 

\begin{remark}[Comparisons to the framework in \cite{balakrishnan2023fundamental}] \label{remark:bala} In our context, the minimax framework introduced by \cite{balakrishnan2023fundamental} is to show that the following holds
\begin{align*}
    \forall g \in \mathcal{T}, \qquad \inf_{\hat{\theta}} \sup_{\mathbb{P}\in \mathcal{P}_\mu(\mathcal{F}_\mu(g))} \mathbb{E}_{\mathbb{P}^{2n}}\left[\left|\hat{\theta} - \theta(\mathbb{P})\right|\right] \asymp \nnhalf + \delta_\mu^\appr 
\end{align*} where $\mathcal{F}_{\mu}(g) = \{\mu_0 \in \mathcal{T}: \|\mu_0 - g\|_2 \le \delta_\mu^\appr\}$. The result is close to our minimax risk $\mathfrak{m}(n, (\delta^\appr_\mu, 0))$ by ignoring the effect of the stochastic error of the function class, i.e., $\delta_\mu^\stoc = 0$, with a slightly stronger instance-level lower bound showing that the same lower bound holds for any fixed $g=\hat{\mu}$, the fixed estimate of $\mu_0$ in their context. Compared with their formalization, we explicitly disentangle the effects of the approximation error and stochastic error because the structure-aware models benefit from mitigating the effect of the stochastic error, as discussed in the introduction. 
\end{remark}

One can directly obtain an upper bound of the above quantity \eqref{eq:sa-risk-mu} from \cref{thm:est1}:  for any $\mathcal{F}_\mu \in \mathcal{H}_n(\bar{\delta}_\mu)$, by the definition of $\mathcal{H}_n(\bar{\delta}_\mu)$, there exists some $\mathcal{G}_\mu$ such that its stochastic error is upper bounded by $\delta_\mu^\stoc$, and $\inf_{g\in \mathcal{G}_\mu} \|f_0 - g\|_2 \le \delta_\mu^\stoc$ for any $f_0 \in \mathcal{F}_\mu$, then applying \cref{thm:est1} using $\mathcal{G}_\mu$ directly gives the upper bound in the following theorem.

\begin{theorem}
\label{thm:sa-mu}
    For any $\bar{\delta}_\mu = (\delta_\mu^\appr, \delta_\mu^\stoc) \in [0, 1]^2$,  the structure-agnostic minimax risk defined in \eqref{eq:sa-risk-mu} with budget $\bar{\delta}_\mu$ satisfies
    \begin{align*}
         \mathfrak{m}(n, \bar{\delta}_\mu) \asymp n^{-1/2} + \delta^\appr_\mu + (\delta^\stoc_\mu)^2.
    \end{align*}
\end{theorem}

For the lower bound, it suffices to identify a single hard class $\mathcal{F}_\mu$ in $\mathcal{H}_n(\bar{\delta}_\mu)$. The next result shows that a sparse linear class, thickened by $L_2$-balls of radius $\delta_\mu^\appr$, is already hard enough.

\begin{proposition}
\label{prop:lb-slm}
    For any $\delta^\appr_\mu \in [0,1]$, $n \ge 18$, and $m, s \in \mathbb{N}^+$ satisfying $m\ge n^{3}$ and $\frac{s[\log(m/s)+1]}{n} \le 1$, there exists a set of orthogonal basis $\phi(X) = (\phi_1(X), \ldots, \phi_m(X))^\top \in \mathbb{R}^m$ with $\mathbb{E}_{X\sim \nu_d}[\phi(X) \phi(X)^\top] = I_m$ such that the set defined as
    \begin{align}
        \label{eq:fmu-text}
        \mathcal{F}_{\mu} := \bigcup_{g\in \mathcal{G}_\mu} \mathcal{B}(g, \delta^\appr_\mu) \qquad \text{with} ~~~ \mathcal{G}_\mu := \left\{g(x) = \sum_{j=1}^m \alpha_j \phi_j(x) \in \tfunc, \|\alpha\|_0 \le s \right\},
    \end{align} and $\mathcal{B}(g, \delta) = \{f \in \tfunc: \|f - g\|_{L_2(\nu_d)} \le \delta\}$ satisfies
    \begin{align}
    \label{eq:lb-1}
         \inf_{\hat{\theta}} \sup_{\mathbb{P}\in \mathcal{P}_\mu(\mathcal{F}_\mu)} \mathbb{E}_{\mathbb{P}^{2n}}\left[|\hat{\theta} - \theta(\mathbb{P})|\right] \ge \frac{1}{211} \left[n^{-1/2} + \delta_\mu^\appr + \frac{s\log(m/s)}{n} \right].
    \end{align} 
\end{proposition}

This proposition gives a lower bound for estimating the linear coefficient when the outcome regression is approximately sparse linear. Compared with prior lower bounds for the well-specified sparse linear model under Gaussian design \citep{cai2017confidence, javanmard2017debiasing}, it explicitly shows the cost of misspecification through $\delta_\mu^\appr$ and works over a uniformly bounded function class without introducing extra logarithmic factors.
Moreover, the class in \eqref{eq:fmu-text} belongs to $\mathcal{H}_n((\delta^\appr_\mu, \delta^\stoc_\mu))$ for any $\delta^\stoc_\mu:=\sqrt{s\log(m/s)/n}$ satisfying $\delta^\stoc_\mu \in [\sqrt{\log(n)/n}, 0.1]$. Hence, \eqref{eq:lb-1} yields the lower bound in \cref{thm:sa-mu} in this range; outside it, the result follows from monotonicity and the parametric term $n^{-1/2}$.
Theorem \ref{thm:sa-mu} shows that when $\pi_0$ cannot be consistently estimated, structure-agnosticity imposes no additional price beyond the worst structure-aware case. The rate
\begin{align*}
n^{-1/2}+\delta_\mu^\appr+(\delta_\mu^\stoc)^2
\end{align*}
is achievable for any black-box ML class with the stated budgets by \cref{thm:est1}, and is already unimprovable for the sparse linear class by \cref{prop:lb-slm}. 

\section{Extension to semiparametric linear functionals}
\label{sec:general}

We extend the upper bound developed for the partial linear model to a much
broader class of semiparametric linear functionals. This framework allows the target functional to depend on a generic and potentially complex conditional mean function and includes the partial linear coefficient in \cref{sec:plm} and common causal estimands, such as the average treatment effect (ATE), as special cases. The message we would like to convey is 
\begin{quote}
\it For a wide range of semi-parametric functional estimation problems, the ``(approximation error) + (stochastic error)$^2$'' error rate is attainable when adopting black-box models for estimating the outcome function, even if it is impossible to estimate the corresponding Riesz representer, a generalization of the auxiliary function, consistently. 
\end{quote}
Thus, in the structure-agnostic regime, the multiplicative rates obtained by standard DML are not fundamental.

\subsection{Problem setup}
Consider the following model of $Y$: 
\begin{align}
\label{eq:model-general}
    Y = m_0(T, X) + \varepsilon \qquad \text{with} \qquad m_0 \in \mathcal{H}_m, ~~~\mathbb{E}[\varepsilon|T, X] = 0,
\end{align}
where $\mathcal{H}_m \subseteq L_2(\nu_z)$ is  a Hilbert space equipped with the inner multiplicative $\langle h_1, h_2\rangle = \mathbb{E}_{\nu_z}[h_1(Z) h_2(Z)]$, and $\nu_z$ is the joint distribution of $Z=(T, X)$.  
Let $\psi: L_2(\nu_z) \times \mathbb{R}^{d+1} \to \mathbb{R}$ be a given function and define the target functional 
\begin{align}
\label{eq:target}
    \theta_0 = \theta(m_0), \quad\mbox{with}\quad \theta(m) = \mathbb{E}_{\nu_z}[\psi(m, Z)].
\end{align}
We assume that $\theta(m)$ is linear and mean-square continuous in $L_2(\nu_z)$, and thus the corresponding Riesz representer is well-defined. This requires the following conditions.

\begin{condition}
\label{cond:linear-functional}
    The following holds for our functional $\theta(\cdot)$:
    \begin{itemize}
    \item[(1)] The functional $\theta(\cdot)$ is linear, i.e., $\theta(h + \alpha \tilde{h}) = \theta(h) + \alpha \theta(\tilde{h})$ for any $h, \tilde{h} \in \mathcal{H}_m$ and $\alpha \in \mathbb{R}$. 
    \item[(2)] The functional $\theta(\cdot)$ is mean-square continuous: there exists a large constant $\const>0$ such that \begin{align*}
        \forall m \in \mathcal{H}_m, \qquad \sqrt{\mathbb{E}_{\nu_z}[|\psi(m, Z)|^2]} \le \const \|m\|_{L_2(\nu_z)}.
    \end{align*}
    \end{itemize}
\end{condition}

Under \cref{cond:linear-functional}, the functional of interests $\theta(\cdot)$ is a continuous linear functional because $|\theta(m)| \le \sqrt{\mathbb{E}_{\nu_z}[|\psi(m, Z)|^2]} \le \const \|m\|_{L_2(\nu_z)}$. By the Riesz representer theorem, there exists a unique Riesz representer of the linear functional $\theta(\cdot)$, denoted as $\repr_0$, in $\mathcal{H}_m$ and satisfies
\begin{align}
\label{eq:a-riesz-rep}
    \forall m\in \mathcal{H}_m \qquad \theta(m) = \mathbb{E}_{\nu_z}[m(Z) \repr_0(Z)]. 
\end{align}
The difficulty is that $b_0(\cdot)$ is unknown and may not be consistently estimated. The partial linear model is a special case with $\mathcal{H}_m = \{\beta \cdot t + \mu(x): \beta \in \mathbb{R}, \mu(x) \in L_2(\nu_x)\}$, $\psi(m) = \beta$, whose corresponding Riesz representer is $\repr_0(t, x) = (t - \mathbb{E}[T|X=x]) / \mathbb{E}[u^2]$.

We observe $2n$ i.i.d. observations of $\{(X_i, T_i, Y_i)\}_{i=1}^{2n}$ from the above model \eqref{eq:model-general} and black-box class $\mathcal{G}_m \subseteq \mathcal{H}_m$ to estimate $m_0$ and subsequently the target parameter $\theta_0 = \theta(m_0)$. We also assume that all the above quantities are uniformly bounded by $\const$.

\begin{condition}[Boundedness]
\label{cond:boundedness}
The following holds for some large constant $\const \ge 1$:
\begin{itemize}
\item[(1)] $\|m_0\|_\infty, \|\repr_0\|_\infty \le \const$, $|\varepsilon| \le \const$.
\item[(2)] $\|g\|_\infty \le \const$ for any $g\in \mathcal{G}_m$.
\item[(3)] $\|\psi(m_0, z)\|_\infty, \|\psi(g, z)\|_\infty \le \const$ for any $g\in \mathcal{G}_m$. 
\end{itemize}
\end{condition}

\begin{example}[Average treatment effect] \label{ex:ate} Let us illustrate the above abstract setup by an example -- the average treatment effect under binary treatment. In this case, the conditional moment model in \eqref{eq:model-general} is instantiated by 
\begin{align}
\label{eq:model-ate}
    Y = \tau_0(X) \cdot T + \mu_0(X) + \varepsilon \qquad \text{with} \qquad T\in \{0, 1\}, ~~~\mathbb{E}[\varepsilon|T, X] = 0,
\end{align} the Hilbert space \begin{align*}
    \mathcal{H}_{m, \ate} := \left\{m(x, t) = t \cdot \tau(x) + \mu(x): \mu,\tau \in L_2(\nu_x) \right\}, \quad \psi_{\mathrm{ATE}}(m; Z) = \tau(X),
\end{align*} 
and the target is $\theta_0 = \mathbb{E}[\tau_0(X)]$. Under overlap and unconfoundedness assumption,  $\theta_0$ coincides with the ATE $\mathbb{E}[Y(1) - Y(0)]$ where $Y(t)$  denotes the potential outcome of $Y$ under treatment $T=t$. The corresponding Riesz representer is
\begin{align}
\label{eq:repr-ate}
    \repr_0(t, x) := t \cdot \frac{1}{\pi_0(x) (1 - \pi_0(x))} - \frac{1}{1 - \pi_0(x)} \qquad \text{with} \qquad \pi_0(x) = \mathbb{E}[T|X=x]
\end{align}
which can be verified from the identity $\theta(m)=\mathbb E_{\nu_z}\{m(Z)\repr_0(Z)\}$. One can easily verify that \cref{cond:linear-functional} and \cref{cond:boundedness} hold when all the components in the model are uniformly bounded and $\pi_0(x) \in [1/C, 1-1/C]$ for some large constant $C>0$. 
\end{example}

\subsection{Method and upper bound result}

We again perform sample splitting and use $\mathcal{D}_2$ to estimate the conditional mean $m_0$  via least squares
\begin{align}
\label{eq:est3-g}
    \hat{g} &\in \argmin_{g\in \mathcal{G}_m} \ninv \sum_{i=n+1}^{2n} (Y_i - g(Z_i))^2 
\end{align}
Using $\mathcal{D}_1$, compute empirical debiasing weights $\hat{a} = (\hat{a}(Z_1),\ldots, \hat{a}(Z_n)) = (\hat{a}_1,\ldots, \hat{a}_n) \in \mathbb{R}^n$ in $\mathcal{D}_1$ by  
\begin{align}
\label{eq:est3-a}
    \hat{a} = \argmin_{\|a\|_\infty \le M} \frac{\lambda}{n} \sum_{i=1}^n a_i^2 + \sup_{f \in \partial \mathcal{G}_m} {\Phi}(a, f)
\end{align} 
where both $M>0$ and $\lambda > 0$ are hyper-parameters to be determined, $\partial \mathcal{G}_m = \{g - \tilde{g}: g,\tilde{g} \in \mathcal{G}_m\}$ and $\Phi(a, f)$ is defined as 
\begin{align}
\label{eq:def-phi}
    \Phi(a, f) = \left|\ninv \sum_{i=1}^n \psi(f, Z_i)- a_i f(Z_i)\right| - \ninv \sum_{i=1}^n \left|f(Z_i)\right|^2.
\end{align} 
Note that when $\psi(\cdot, Z)$ is differentiable in $f$, one can also use gradient descent ascent to solve the above minimax objective, akin to the computation for the partial linear model. Our final SADE estimate of $\theta_0 = \mathbb{E}[\psi(m_0, Z)]$ is the debiased estimator on $\mathcal{D}_1$:
\begin{align}
\label{eq:est3-theta}
    \hat{\theta} = \ninv \sum_{i=1}^n \psi(\hat{g}, Z_i) + (Y_i - \hat{g}(Z_i)) \cdot \hat{a}_i.
\end{align}

\begin{remark}[Simplification when $\psi$ is odd in $f$] In many cases when $\psi(f, Z)$ is odd in $f$, i.e., $\psi(f, Z)= -\psi(-f, Z)$ for any $f \in \mathcal{H}_m$, we have $\sup_{f \in \partial \mathcal{G}_m} \Phi(a, f) = \sup_{f \in \partial \mathcal{G}_m} \Phi_0(a, f)$ with 
\begin{align*}
\Phi_0(a, f) = \ninv \sum_{i=1}^n \psi(f, Z_i)- a_i f(Z_i) - \ninv \sum_{i=1}^n \left|f(Z_i)\right|^2.
\end{align*} This is more convenient in computation, since the resulting objective function continues to be efficiently optimized by stochastic gradient descent ascent (with unbiased loss and gradient) when $n$ is large. 
\end{remark}

For general linear functionals, the stochastic error must control both the
localized complexity of $\partial\mathcal G_m$ and that of its image under
$\psi$. Define $(\partial \mathcal{G}_m)^\psi =  \{\psi(g - \tilde{g}, Z): g, \tilde{g}\in \mathcal{G}_m\}$ and 
\begin{align}
\label{eq:def-stoc-general}
    \delta_m^\stoc := \inf \left\{\delta_s \ge \sqrt{\frac{\log(n)}{n}}: \sup_{\delta \ge \delta_s} \frac{\mathsf{R}_{n, \nu_z}(\delta; \partial \mathcal{G}_m) \lor \mathsf{R}_{n, \nu_z}(\delta; (\partial \mathcal{G}_m)^\psi)}{\const \cdot \delta} \le \delta_s \right\}.
\end{align} The definition of the approximation is similar to that of the partial linear model:
\begin{align}
\label{eq:def-appr-general}
    \delta^\appr_m := \inf_{g\in \mathcal{G}_m} \|m_0 - g\|_{L_2(\nu_z)}.
\end{align}

The following theorem provides an oracle-type inequality for the estimator $\hat{\theta}$ \eqref{eq:est3-theta} we constructed. 

\begin{theorem}
\label{thm:main}
    Under \cref{cond:linear-functional} and \cref{cond:boundedness}, there exists a constant $\tilde{C} = \poly(\const)$ such that with probability at least $1-5e^{-t} - 2\exp(-n(\delta_m^\stoc)^2)$, if $M\ge \const$, our estimator $\hat{\theta}$ in \eqref{eq:est3-theta} satisfies
    \begin{align}\label{eq:rate-main}
        |\hat{\theta} - \theta_0| \le \tilde{C} \cdot M \left(\sqrt{\frac{t}{n}} + \delta_m^\appr + (\delta_m^\stoc)^2 + \lambda \right)
    \end{align} with the notations $\delta_m^\stoc$ and $\delta_m^\appr$ defined in \eqref{eq:def-stoc-general} and \eqref{eq:def-appr-general}. 
\end{theorem}

\cref{thm:main} extends the upper bound result in \cref{thm:est1} for generic linear functional estimation, and many implications thus follow similarly. We make a summary below in the regime where it is impossible to consistently estimate the Riesz representer. One can also derive a similar convergence result akin to \cref{prop:nn-plm} for the partial linear model; see that for the example of ATE \cref{ex:ate} in \cref{sec:ex-ate}. 
\begin{itemize}
    \item Compared with the DML rate $n^{-1/2} + \delta_m^\appr + \delta_m^\stoc$, our method removes the first-order dependency on the stochastic error $\delta_m^\stoc$ by the empirical debiasing weights in \eqref{eq:est3-a} targeted at $\partial \mathcal G_m$. 
    \item We should also do under-smoothing $\delta_m^\appr \asymp (\delta_m^\stoc)^2$ in model selection for $\mathcal{G}_m$ for optimal estimation of $\theta_0$. In this case, root-n consistency is still possible when there exists a black-box model class satisfying $\delta_m^\appr \lesssim n^{-1/2} \text{ and } \delta_m^\stoc \lesssim n^{-1/4}$.
\end{itemize}

\appendix
\appendixrefs
\vspace{15pt}
\noindent \textbf{{\LARGE Appendix}}


\section{Further results for partial linear model}
\label{sec:further}

In this section, we provide detailed theoretical analyses for the claims we argue in the main text for completeness. In \cref{sec:plm-full} we introduce the structure-agnostic minimax risk when both $\mu_0$ and $\pi_0$ can be estimated by machine learning models and present the tightest upper and lower bounds we can obtain now. 

\subsection{Full structure-agnostic minimax risk}
\label{sec:plm-full}

Recall the definition of $\mathcal{H}_n(\bar{\delta})$ in \eqref{eq:h-n-delta}. Given two specific function classes $\mathcal{F}_\mu$ and $\mathcal{F}_\pi$, we define the family of distributions as 
\begin{align}
\label{eq:dgp-plm-full}
\begin{split}
    \mathcal{P}(\mathcal{F}_{\mu}, \mathcal{F}_\pi) = \Big\{ &(X, T, Y) \sim \mathbb{P}_{\beta, v, \mu,\pi}~\text{where}~ \mu \in \mathcal{F}_{\mu}, \pi\in \mathcal{F}_\pi, |\beta| \lor |v|\le 3, |T|\lor |Y| \le 3:\\
    &~~~~X\sim \nu_d, ~~ T = \pi(X) + u ~\text{with}~ \mathbb{E}[u|X]=0,\mathbb{E}[u^2|X] \equiv v \ge 1/3\\
    &~~~~Y = T \cdot \beta + \mu(X) + \varepsilon ~\text{with}~ \mathbb{E}[\varepsilon_Y|X, T]=0\Big\},
\end{split}
\end{align}

Now we define the structure-agnostic minimax risk with approximation and stochastic error budgets $(\bar{\delta}_\mu, \bar{\delta}_\pi)$ for $(\mu, \pi)$ as 
\begin{align}
\label{eq:sa-risk}
    \mathfrak{m}(n, \bar{\delta}_\mu, \bar{\delta}_\pi) := \sup_{\substack{\mathcal{F}_\mu \in \mathcal{H}_n(\bar{\delta}_\mu)\\ \mathcal{F}_\pi \in \mathcal{H}_n(\bar{\delta}_\pi)}} \underbrace{\inf_{\hat{\theta}} \sup_{\mathbb{P}\in \mathcal{P}(\mathcal{F}_\mu, \mathcal{F}_\pi)} \mathbb{E}_{\mathbb{P}^{2n}}\left[\left|\hat{\theta} - \theta(\mathbb{P})\right|\right]}_{\text{standard minimax risk over given class $(\mathcal{F}_\mu, \mathcal{F}_\pi)$}},
\end{align} where the randomness in $\mathbb{E}_{\mathbb{P}^{2n}}$ is $2n$ i.i.d. observations from $\mathbb{P}$, and $\hat{\theta}$ is the function of the $2n$ observations and potential black-box machine learning models $\mathcal{G}_\mu, \mathcal{G}_\pi$ that are tailored to $\mathcal{F}_\mu$ and $\mathcal{F}_\pi$, respectively. 

The upper bound of the minimax risk \eqref{eq:sa-risk} is obtained by combining our estimator and the standard double machine learning estimator because the algorithm itself has the knowledge of the budgets $\bar{\delta}_\mu$ and $\bar{\delta}_\pi$. The lower bounds in \eqref{eq:sa-risk} are obtained by choosing the hard function classes $\mathcal{F}$ as (1) $L_\infty$ balls with radius $\delta$ and (2) sparse linear models. 
\begin{theorem}
\label{thm:sa}
    For any $\bar{\delta}_\mu, \bar{\delta}_\pi \in [0, 1]^2$,  the structure-agnostic minimax risk defined in \eqref{eq:sa-risk} with budgets $\bar{\delta}_\mu, \bar{\delta}_\pi$ satisfies
    \begin{align*}
         \mathfrak{m}(n, \bar{\delta}_\mu, \bar{\delta}_\pi) \lesssim n^{-1/2} + \underbrace{\left[\delta^\appr_\mu + (\delta^\stoc_\mu)^2\right]}_{\text{Rate by our method}} \wedge \underbrace{\left[(\delta^\appr_\mu + \delta^\stoc_\mu) \cdot (\delta^\appr_\pi + \delta^\stoc_\pi)\right]}_{\text{Rate by standard DML}},
    \end{align*} and 
    \begin{align*}
        \mathfrak{m}(n, \bar{\delta}_\mu, \bar{\delta}_\pi) \gtrsim n^{-1/2} + \left[\delta^\appr_\mu + (\delta^\stoc_\mu)^2\right] \cdot \left[\delta^\appr_\pi + (\delta^\stoc_\pi)^2\right] + \left[\delta^\stoc_\mu \land \delta^\stoc_\pi\right]^2.
    \end{align*}
\end{theorem}

\section{Result for average treatment effect estimation}
\label{sec:ex-ate}

We instantiate the ML model class $\mathcal{G}_m$ using structured ReLU NNs:
\begin{align}
\label{eq:gm-ate}
    \mathcal{G}_{m,\ate} = \left\{m(t, x) = t \cdot g_1(x) + g_2(x): g_1, g_2 \in \mathcal{F}_{\mathtt{nn}}(d, L, N, B)\right\}.
\end{align}

\begin{condition} \label{cond:nn-ate} Under  model \eqref{eq:model-ate} with $\pi_0(x)=\mathbb{E}[T|X=x]$, there exists a constant $c_2 > 1$ such that the following holds:
\begin{itemize}
\item[(1)] (Function complexity) $\mu_0 \in \mathcal{F}_{\mathtt{HCM}}(d, l, \mathcal{O}_\mu, C)$, $\tau_0 \in \mathcal{F}_{\mathtt{HCM}}(d, l, \mathcal{O}_\tau, C)$ with $\max_{(\beta, t) \in \mathcal{O}_\mu \cup \mathcal{O}_\tau} (\beta\lor t) \lor l \lor C \lor d \le c_2$. Let $\gamma_\mu^\star = \min_{(\beta, t)\in \mathcal{O}_\mu} (\beta/t)$, $\gamma_\tau^\star = \min_{(\beta, t)\in \mathcal{O}_\tau} (\beta/t)$
\item[(2)] (Boundedness) All the components in the model are bounded: $\|X\|_\infty \le 1$, $|\varepsilon| \lor |u| \le c_2$, $\|\mu_0\|_\infty \lor \|\tau_0\|_\infty$.
\item[(3)] (Non-overlap) The propensity score is bounded away from $0$ and $1$: $\|\pi_0^{-1}\|_\infty \lor \|(1 - \pi_0)^{-1}\|_\infty \le c_2$. 
\item[(4)] (Machine learning class and hyper-parameter $M$) We adopt $\mathcal{G}_m \gets \mathcal{G}_{m,\ate}$ as in \eqref{eq:gm-ate} with $B=c_2$, and pick $M=2c_2$. 
\end{itemize}
\end{condition}

\begin{proposition}\label{prop:nn-ate}
Under \cref{cond:nn-ate}, our estimator $\hat{\theta}$ in \eqref{eq:est3-theta} with NN hyperparameter satisfying $N\land L\ge \tilde{C}_2(1+\log(n))$  and $NL \asymp n^{\frac{1}{2[(\gamma_\mu^\star\land \gamma_\tau^\star)+1]}} \log^{\frac{2(\gamma^\star_\mu\land \gamma_\tau^\star)-1}{(\gamma^\star_\mu \land \gamma_\tau^\star)+1}}(n)$ satisfies, with probability at least $1-n^{-10}$:
\begin{align}
\label{eq:rate-nn-ate}
    |\hat{\theta} - \theta_0| \le \tilde{C}_2 \left[\lambda + \sqrt{\frac{\log(n)}{n}} +  \left(\frac{\log^6(n)}{n}\right)^{\frac{(\gamma_\mu^\star \land \gamma_\tau^\star)}{(\gamma_\mu^\star \land \gamma_\tau^\star)+1}}\right].
\end{align} for some constant $\tilde{C}_2$ dependent on $c_2$. 
\end{proposition}

The optimal complexity choice is obtained by under-smoothing in a similar manner to the partial linear model. The effective smoothness $(\gamma_\mu^\star \land \gamma_\tau^\star)$ is governed by the harder of the two regression components, $\mu_0$ and $\tau_0$, and the resulting rate follows from balancing approximation error against the \emph{square} of stochastic error rather than against stochastic error itself.

\bibliographystyle{apalike2}
\bibliography{statbb.bib}

\clearpage
\setcounter{section}{0}
\setcounter{subsection}{0}
\setcounter{equation}{0}
\setcounter{theorem}{0}

\renewcommand{\thesection}{S\arabic{section}}
\renewcommand{\thesubsection}{\thesection.\arabic{subsection}}
\renewcommand{\theequation}{\thesection.\arabic{equation}}
\renewcommand{\thetheorem}{\thesection.\arabic{theorem}}

}
\end{document}